\documentclass[12pt]{article}
\usepackage{mathrsfs}
\usepackage{amsmath}
\usepackage{amssymb}
\usepackage{amsfonts}
\usepackage{wasysym}
\usepackage{verbatim}
\usepackage[amsmath,thmmarks]{ntheorem}

\def\F{\mathscr{F} }

\def\R{\mathbb{R}}
\def\Z{\mathbb{Z}}
\def\N{\mathbb{N}}
\def\T{\mathbb{T}}

\def\ba{\begin{array}}
\def\ea{\end{array}}
\def\be{\begin{enumerate}}
\def\ee{\end{enumerate}}
\def\bi{\begin{itemize}}
\def\ei{\end{itemize}}
\def\bc{\begin{cases}}
\def\ec{\end{cases}}
\def\beq{\begin{equation}}
\def\eeq{\end{equation}}
\def\beqs{\begin{equation*}}
\def\eeqs{\end{equation*}}
\def\beqa{\begin{eqnarray}}
\def\eeqa{\end{eqnarray}}
\def\beqas{\begin{eqnarray*}}
\def\eeqas{\end{eqnarray*}}
\def\bmul{\begin{multline}}
\def\emul{\end{multline}}
\def\bmuls{\begin{multline*}}
\def\emuls{\end{multline*}}
\def\bg{\begin{gather}}
\def\eg{\end{gather}}
\def\bgs{\begin{gather*}}
\def\egs{\end{gather*}}

\theorembodyfont{\normalfont\sl}
\newtheorem{thm}{Theorem}[section]
\newtheorem{cor}[thm]{Corollary}
\newtheorem{lem}[thm]{Lemma}
\newtheorem{prop}[thm]{Proposition}

\numberwithin{equation}{section}

\begin{document}

%\baselineskip=24pt \thispagestyle{empty}
%
%\baselineskip=24pt
%
%\mbox{}
%\bigskip
%
%\centerline{\LARGE\bf  Global Attractor for Weakly Damped}
%
%\bigskip
%
%\centerline{\LARGE\bf  Forced KdV Equation in}
%
%\bigskip
%
%\centerline{\LARGE\bf   Low Regularity}
%
%
%\bigskip
%
%
%\centerline{Yifei WU}
%
%\centerline{Department of Mathematics, South China University
%            of  Technology,}
%
%\centerline{
%            Guangzhou, Guangdong 510640, P. R. China}
%
%\bigskip
%\bigskip
%
%\leftskip 1 true cm \rightskip 1 true cm
%
% \noindent {\bf Abstract} \
%In this paper we consider the long time behavior of the weakly
%damped, forced Korteweg-de Vries equation in the Sololev spaces of
%the negative indices in the periodic case. We prove that the
%solutions are uniformly bounded in $\dot{H}^s(\T)$ for $s>-\dfrac{1}{2}$.
%Moreover, we show that the solution-map possesses a global attractor
%in $\dot{H}^s(\T)$, which is a compact set in $H^{s+3}(\T)$.
%
%\bigskip
%
%\noindent {\bf Keywords:}\quad Korteweg-de Vries equation,
%attractor, asymptotic smoothing effect, Bourgain space, I-method.
%
%\bigskip
%
%\noindent {\bf  MR(2000) Subject Classification:}\quad 35Q53, 35B41
%\bigskip
%
%\noindent {\bf Short Title:} Attractor for the KdV Equation
%
%
%\bigskip
%\noindent {\bf Corresponding author:} Yifei WU
%
%
%\bigskip
%\noindent {\bf Email Address:}  yerfmath@yahoo.cn
%
%\leftskip 0 true cm \rightskip 0 true cm

\newpage

 \baselineskip=24pt

\setcounter{page}{1}

\title{
\baselineskip=24pt
 \bf Global Attractor for Weakly Damped\\[.5ex]
 \bf Forced KdV Equation in\\[.5ex]
 \bf Low Regularity on $\T$\footnote{This work is
 supported by National Natural Science Foundation of China
 under grant numbers 10471047 and 10771074.}}

\bigskip

\author{
\baselineskip=24pt
 \normalsize Yongsheng Li \quad and \quad
             Yifei Wu\footnote{Email: yshli@scut.edu.cn (Y. S. Li) and
             yerfmath@yahoo.cn (Y. F. Wu)}\\[1ex]
 \normalsize Department of Mathematics,
             South China University of Technology, \\[1ex]
 \normalsize Guangzhou, Guangdong 510640, P. R. China
 }

\bigskip

\date{}

\bigskip

\bigskip

\leftskip 1cm \rightskip 1cm

\maketitle

\noindent
 {\small{\bf Abstract}\quad
 \baselineskip=24pt
In this paper we consider the long time behavior of the weakly
damped, forced Korteweg-de Vries equation in the Sololev spaces of
the negative indices in the periodic case. We prove that the
solutions are uniformly bounded in $\dot{H}^s(\T)$ for
$s>-\dfrac{1}{2}$. Moreover, we show that the solution-map possesses
a global attractor in $\dot{H}^s(\T)$ for $s>-\dfrac{1}{2}$, which
is a compact set in $H^{s+3}(\T)$.}
\bigskip

\noindent
 {\small{\bf Keywords:} Korteweg-de Vries equation,
attractor, asymptotic smoothing effect, Bourgain space, I-method.}

\bigskip

\noindent
 {\small{\bf  MR(2000) Subject Classification:}\quad 35Q53,
35B41}

\bigskip

\bigskip

\leftskip 0cm \rightskip 0cm

\normalsize

\baselineskip=24pt

\section{Introduction}

In this paper we study the long time behavior of the weakly damped,
forced Korteweg-de Vries (KdV) equations
\renewcommand{\arraystretch}{2}
\begin{eqnarray}\label{KdV}
  &&\partial_{t}u+\partial_{x}^{3}u+\dfrac{1}{2}\partial_xu^2+\gamma u
  =f,\qquad x\in \T,\,t\in \R^+,\\
  &&u(x,0)=u_0(x)\in \dot{H}^s(\T),\label{initial data}
\end{eqnarray}
where the unknown function $u$ is real-valued, the external force
$f$ is time-independent and belongs to $\dot{H}^s(\T)$, and the
damping parameter $\gamma>0$. Here $\T=\R/2\pi \Z$ and the space
$$
\dot{H}^s(\T)=\left\{f\in H^s(\T): f\mbox{ is zero mean, } i.e.
\displaystyle \int_\T f(x) \,dx =0\right\}.
$$

When $\gamma=0$, $f=0$, it is known as KdV equation, which has been
widely studied. The local well-posedness is known in $H^s(\T)$ for
$s\geq -1/2$ by Kenig, Ponce and Vega \cite{KPV1} by using an
improved Bourgain argument. In \cite{CKSTT2}, Colliander, Keel,
Staffilani, Takaoka and Tao introduced the well-known I-method to
extended the local result to a global one when $s\geq-1/2$. Further,
there are examples given in \cite{CCT} to show that the KdV equation
is ill-posed for $s<-1/2$. On the other hand, the global attractor
for the weakly damped, forced KdV equation was studied by many
authors, see \cite{Gh2}--\cite{Goubet}, \cite{GR}, \cite{GH2}
\cite{L}, \cite{MoRo}--\cite{Ts}, $ect.$. Particularly, in
\cite{Goubet}, Goubet proved that the existence of the global
attractor in $\dot{L}^2(\T)$, which is a compact set in $H^3(\T)$.
He showed these by the energy equations method (J.\ M.\ Ball's
argument) and Bourgain's argument. In \cite{Ts}, Tsugawa consider
(\ref{KdV}) with $f\in \dot{L}^2(\T)$, and prove the existence of
the global attractor in $\dot{H}^s(\T)$ for $s>-\dfrac{3}{8}$, which
is compact in $H^3(\T)$. However, it seems that the method used in
\cite{Ts} depends strongly on the assumption $f\in L^2(\T)$ to
obtain the uniform boundness of the solution, moreover, it has a gap
about the indix $s$ compared to the local well-posedness for
$s\geq-\dfrac{1}{2}$. Meanwhile, the argument in \cite{Ts} seems not
suitably to the non-compact case.

In this paper, we consider (\ref{KdV}) and prove that the global
attractor exists in $\dot{H}^s(\T)$ for $s>-\dfrac{1}{2}$, and is
compact in $H^{s+3}(\T)$, which was left open in \cite{Goubet} and
\cite{Ts}. Our main approach is combining the method in
\cite{Goubet} and the idea of I-method. We first use the I-method to
prove the uniform boundness of the solutions in Section 3. We remark
that the argument used here should be somewhat different from
\cite{CKSTT2}, since what we need here is the uniform boundness of
the solutions. However, it can be also succeeded in applying here,
based on the brief that the weak damping shall prevent the
unlimitedly increasing of the the energy in $H^s(\T)$. Compared to
\cite{Ts}, the strong condition $f\in \dot{L}^2(\T)$ is relaxed to
general: $f\in \dot{H}^s(\T)$, and the index is lowered to be
optimal by this argument. As a consequence, we obtain the existence
of the bounded absorbing sets in $H^s(\T)$. Furthermore, we note
that it lacks the conservation law in $H^s(\T)$ when $s<0$ for the
KdV equation, so it is not sufficient to use the energy equation
method to show the asymptotic compactness in $H^s(\T)$. To overcome
this problem, we pursue the asymptotic smoothing effect of the
flow-map at first, and then prove the asymptotic compactness by
using the asymptotic smoothing effect. More precisely, we split the
solution into two parts, and prove that one part is decay to zero
and the other is regular (bounded in $H^{s+3}(\T)$), which is shown
in Section 4. As we know, we need some infinitesimals (we denote
them by $\epsilon(N)$ below, see Section 4 for the details), which
indicate that the energies from nonlinearities can be controlled to
ensure the decay of one part and uniformly bounded in a higher
regular space of the other. However, it is not easy to give these
$\epsilon(N)$ in this situation when the solutions are
distributions. To cope with the difficulties, some decomposition
tricks to treat the nonlinearities and special multilinear estimates
included some commutators will be used throughout the article.
Furthermore, some arguments in \cite{M} and \cite{Ts2} are employed
in this section. By employing this asymptotic smoothing effect, one
may prove directly that the solution map is asymptotic compact in
$H^{l}(\T)$ for any $l<s+3$. At last, we further show the
compactness of the global attractor in $H^{s+3}(\T)$. We leave the
existence of the global attractor in critical space
$\dot{H}^{-\frac{1}{2}}(\T)$ open in this paper.

{\bf Notations.} We use $A\lesssim B$ or $B\gtrsim A$ to denote
the statement that $A\leq CB$ for some large constant $C$ which
may vary from line to line, and may depend on the data such as
$\gamma,f,u_0$ and the index $s$ unless otherwise mentioned. When it is
necessary, we will write the constants by $C_1,C_2,\cdots$, or
$K_1,K_2,\cdots$ to see the dependency relationship. We use $A\ll
B$ to denote the statement $A\leq C^{-1}B$, and use $A\sim B$ to
mean $A\lesssim B\lesssim A$. The notation $a+$ denotes
$a+\epsilon$ for any small $\epsilon$, and $a-$ for $a-\epsilon$.
$\langle\cdot\rangle=(1+|\cdot|^2)^{1/2}$,
$J_x^\alpha=(1-\partial^2_x)^{\alpha/2}$ and
$D_x^\alpha=(-\partial^2_x)^{\alpha/2}$. We use
$\|f\|_{L^p_xL^q_t}$ to denote the mixed norm
$\Big(\displaystyle\int\|f(x,\cdot)\|_{L^q}^p\
dx\Big)^{\frac{1}{p}}$. Moreover, we denote $\F_x$ to be the
Fourier transform corresponding to the variable $x$. Further we
define the Fourier spectral projector $P_N$, $Q_N$, $P_{\ll N}$
respectively as
$$
\begin{array}{c}
P_Nf(x)=\displaystyle\int_{|\xi|\leq N}\!
e^{ix\xi}\hat{f}(\xi)\,d\xi,\,\,
Q_Nf(x)=\displaystyle\int_{|\xi|\geq N}\! e^{ix\xi}\hat{f}(\xi)\,d\xi,\,\,
P_{\ll N}f(x)=\displaystyle\int_{|\xi|\ll N}\!
e^{ix\xi}\hat{f}(\xi)\,d\xi.
\end{array}
$$

We state our main result, as the end of this section.
\begin{thm}
Let $-\dfrac{1}{2}<s<0$ and $\gamma>0$, $f\in \dot{H}^s(\T)$. Then
(\ref{KdV}) (\ref{initial data}) possess a global attractor in
$\dot{H}^{s}(\T)$, which is compact in $H^{s+3}(\T)$.
\end{thm}

 \vspace{0.3cm}
\section{Functional Spaces and Preliminary Estimates}

We first introduce some notations and definitions. We define
$(d\xi)_\T$ to be the normalized counting measure on $\T$ such that
$$
\displaystyle\int a(\xi)\,(d\xi)_\T=\frac{1}{2\pi}\sum\limits_{\xi\in
\Z} a(\xi).
$$
At the following, we always adopt the notation
$$\displaystyle\int_\ast a(\xi,\xi_1,\xi_2,\tau,\tau_1,\tau_2)
=\int_{\stackrel{\xi\xi_1\xi_2\neq 0,\xi_1+\xi_2=\xi,} {
\tau_1+\tau_2=\tau}}a(\xi,\xi_1,\xi_2,\tau,\tau_1,\tau_2)\,(d\xi_1)_\T
(d\xi_2)_\T \,d\tau_1 d\tau_2.
$$

Define the Fourier transform of a function $f$ on $\T$ by
$$
\hat{f}(\xi)=\displaystyle\int_\T e^{-2\pi i  x\xi}f(x)\,dx,
$$
and thus the Fourier inversion formula
$$
f(x)=\displaystyle\int e^{2\pi i  x\xi} \hat{f}(\xi)\,(d\xi)_\T.
$$

For $s,b\in \R $, we define the Bourgain space $X_{s,b}$ to be the
completion of the Schwartz class under the norm
\begin{equation}
\|f\|_{X_{s,b}}\equiv \left(\int\!\!\!\!\int
\langle\xi\rangle^{2s}\langle\tau-\xi^3\rangle^{2b}|\hat{f}(\xi,\tau)|^2
\,(d\xi)_\T\, d\tau\right)^{\frac{1}{2}}.\label{X}
\end{equation}
For an interval $\Omega$, we define $X_{s,b}^{\Omega}$ to be the
restriction of $X_{s,b}$ on $\T\times\Omega$ with the norm \beq
\|f\|_{X_{s,b}^\Omega}=\inf\{\|F\|_{X_{s,b}}:F|_{t\in\Omega}=f|_{t\in\Omega}\}.
\label{X1} \eeq When $\Omega=[-\delta,\delta]$, we write
$X_{s,b}^\Omega$ as $X_{s,b}^\delta$. By the Functional analysis on
Hilbert space we see that, for every $f\in X_{s,b}^{\Omega}$, there
exists an extension $\tilde{f}\in X_{s,b}$ such that $\tilde{f}=f$
on $\Omega$ and
$\|f\|_{X_{s,b}^{\Omega}}=\big\|\tilde{f}\big\|_{X_{s,b}}$.

Further, we define the locally well-posed working spaces $Y^s$ via
the norms
$$
\|f\|_{Y^s}= \|f\|_{X_{s,\frac{1}{2}}}
+\|\langle\xi\rangle^s\hat{f}(\xi,\tau)\|_{L^2((d\xi)_\T)L^1(d\tau)},
$$
and the companion spaces $Z^s$:
$$
\|f\|_{Z^s}= \|f\|_{X_{s,-\frac{1}{2}}}
+\left\|\frac{\langle\xi\rangle^s}{\langle\tau-\xi^3\rangle}
\hat{f}(\xi,\tau)\right\|_{L^2((d\xi)_\T)L^1(d\tau)}.
$$
Then it is easy to see that $ Y^s\subset C(\R;H^s(\T))\cap
X_{s,\frac{1}{2}}$.

 Let $s<0$ and $N\gg 1$ be fixed,  the Fourier multiplier
operator $I_{N,s}$ is defined as \beq
\widehat{I_{N,s}f}(\xi)=m_{N,s}(\xi)\hat{f}(\xi),\label{I} \eeq
where the multiplier $m_{N,s}(\xi)$ is a smooth, monotone function
satisfying $0<m_{N,s}(\xi)\leq 1$ and
 \beq m_{N,s}(\xi)=\Biggl\{
\begin{array}{ll}
1,&|\xi|\leq N,\\
N^{-s}|\xi|^s,&|\xi|>2N.\label{m}
\end{array}
\eeq Sometimes we denote $I_{N,s}$ and $m_{N,s}$ by $I$ and $m$
respectively for short if there is no confusion.

It is obvious that the operator $I_{N,s}$ maps $H^s(\T)$ into
$L^2(\T)$ with equivalent norms for any $s<0$. More precisely, there
exists some positive constant $C$ such that \beq
    C^{-1}\|f\|_{H^s}\leq \|I_{N,s}f\|_{L^2}
\leq
    CN^{-s}\|f\|_{H^s}.\label{II}
\eeq

Next, we give some estimates which will be used in the following
sections.
\begin{lem}\!\!\cite{B}.
Let $f\in X_{0,\frac{1}{3}}$, then
\begin{equation}
\left\|f\right\|_{L^4_{xt}}
    \lesssim
    \|f\|_{X_{0,\frac{1}{3}}}.\label{L5.1}
\end{equation}
\end{lem}

\begin{lem}
Let $\varphi\in H^{\frac{1}{2}-}(\R)$, $f\in
X_{s,\frac{1}{2}}$ for any $s\in \R$, then
\begin{equation}\label{5.8}
\|\varphi f\|_{X_{s,\frac{1}{2}-}} \leq
    C(\varphi)\|f\|_{X_{s,\frac{1}{2}}}.
\end{equation}
\end{lem}
\begin{proof}
Denote $U(t)\equiv e^{-\partial_x^3 t}$, then by H\"{o}lder's and
Sobolev's inequalities, we have
$$
\begin{array}{ll}
\|\varphi(t)U(-t) f(x,t)\|_{H^{\frac{1}{2}-}_t} &\!\lesssim
\left\|J_t^{\frac{1}{2}-}\!(\varphi)\> U(-t) f\right\|_{L^2_t}+
\left\|\varphi\> J_t^{\frac{1}{2}-}\!(U(-t)
f)\right\|_{L^2_t}\\
&\!\lesssim \left\|J_t^{\frac{1}{2}-}\!(\varphi)\right\|_{L^p_t}
\left\|U(-t) f\right\|_{L^{p'}_t}
+\left\|\varphi\right\|_{L^{p'}_t}\left\|J_t^{\frac{1}{2}-}U(-t) f\right\|_{L^p_t}\\
&\!\lesssim \|\varphi\|_{H^{\frac{1}{2}-}}
\left\|J_t^{\frac{1}{2}}U(-t) f\right\|_{L^2_t},
\end{array}
$$
where we set $p=2+$ such that $H_t^{0+}(\R)\hookrightarrow
L^p_t(\R) $, and $\dfrac{1}{p}+\dfrac{1}{p'}=\dfrac{1}{2}$. Thus
$H_t^{\frac{1}{2}-}(\R)\hookrightarrow L^{p'}_t(\R)$. Then the
result follows by noting $ \|\varphi f\|_{X_{s,b}}
=\|\varphi(t)U(-t) f(x,t)\|_{H^s_x H^{b}_t} $ for any $s, b\in\R$.
\hfill$\Box$
\end{proof}

{\it \noindent Remark.} As a candidate, we keep in mind that
$\varphi=\chi_{[0,\delta]}(t)e^{-\delta t}$ for some $\delta>0$. A
more precise version of Lemma 2.2 is, for any $0\leq b<b'\leq
\dfrac{1}{2}$,
$$
\|\varphi f\|_{X_{s,b}} \leq
    C(\varphi)\|f\|_{X_{s,b'}}.
$$

We denote $\psi(t)$ to be an even smooth bump function of the
interval $[-1,1]$.
\begin{lem}
For any $u,v\in X_{-\frac{1}{2},\frac{1}{2}}^\delta,$ such that $u,v$
are zero $x-$mean for all $t$, then
\begin{equation}\label{L21}
 \|\psi(t/\delta)\>\partial_x(uv)\|_{Z^{-\frac{1}{2}}}\lesssim
 \delta^{\frac{1}{3}-}\|u\|_{X_{-\frac{1}{2},\frac{1}{2}}^\delta}
 \|v\|_{X_{-\frac{1}{2},\frac{1}{2}}^\delta}.
\end{equation}
\end{lem}
\begin{proof}
It is actually a more precise result from the proof of Proposition 3
in \cite{CKSTT2}, where we shall combine the estimate
\begin{equation}\label{delta}
\|\psi(t/\delta)f\|_{X_{s,b'}} \lesssim
\delta^{b-b'}\|f\|_{X_{s,b}},
\end{equation}
for any $s\in \R,0\leq b'\leq b<\dfrac{1}{2}, f\in X_{s,b}$.
\hfill$\Box$
\end{proof}

{\it \noindent Remark.} Another version from the proof of
Proposition 3 in \cite{CKSTT2} is
\begin{equation}\label{L21''}
 \|\varphi\psi(t/\delta)\>\partial_x(uv)\|_{X_{-\frac{1}{2},-\frac{1}{2}}}
 \lesssim
 \delta^{\frac{1}{3}-}\|u\|_{X_{-\frac{1}{2},\frac{1}{2}}^\delta}
 \|v\|_{X_{-\frac{1}{2},\frac{1}{2}}^\delta}
\end{equation}
for any $\varphi\in H^{\frac{1}{2}-}(\R)$, which
follows by the estimate (\ref{5.8}).

By (\ref{L21''}) and the relation $|\xi|\leq |\xi_1|+|\xi_2|$ when
$\xi=\xi_1+\xi_2$, it leads to
\begin{equation}\label{L21'''}
 \|\varphi\psi(t/\delta)\>\partial_x(uv)\|_{Z^s}
 \lesssim
 \|u\|_{X_{s,\frac{1}{2}}^\delta}\|v\|_{X_{-\frac{1}{2},\frac{1}{2}}^\delta}
 +\|u\|_{X_{-\frac{1}{2},\frac{1}{2}}^\delta}\|v\|_{X_{s,\frac{1}{2}}^\delta},
\end{equation}
and thus
\begin{equation}\label{L21'}
 \|\varphi\psi(t/\delta)\>\partial_x(uv)\|_{X_{s,-\frac{1}{2}}}
 \lesssim
 \|u\|_{X_{s,\frac{1}{2}}^\delta}\|v\|_{X_{-\frac{1}{2},\frac{1}{2}}^\delta}
 +\|u\|_{X_{-\frac{1}{2},\frac{1}{2}}^\delta}\|v\|_{X_{s,\frac{1}{2}}^\delta},
\end{equation}
for any $s\geq-\dfrac{1}{2}$, where we have dropped the factor
$\delta^{\frac{1}{3}-}$ since it is smaller than $1$.

Recall that the commutator $[A,B]=AB-BA$, we have the following
result.
\begin{lem}
For any $s\in \R$, if the functions $u,v$ are zero $x-$mean for all
$t$ and $u=P_{ N}u,v=Q_N v$ for some $N>0$, then
\begin{equation}
\left\|\partial_x[J_x^s,u]v\right\|_{X_{\frac{1}{2},-\frac{1}{2}}}
    \lesssim
    \|u\|_{X_{\frac{1}{2},\frac{1}{2}}}\>
    \|v\|_{X_{s-\frac{1}{2},\frac{1}{2}}}.\label{L5.4}
\end{equation}
\end{lem}
\begin{proof}
By duality and Plancherel's identity, it suffices to show
\begin{eqnarray*}
    &&\displaystyle\int_{\ast}
    \xi\langle\xi\rangle^{\frac{1}{2}}
    \langle\xi_1\rangle^{-\frac{1}{2}}
    \langle\xi_2\rangle^{-s+\frac{1}{2}}
    (\langle\xi\rangle^{s}-\langle\xi_2\rangle^{s})\>
    \hat{f}(\xi,\tau)\>\hat{g}(\xi_1,\tau_1)\>\hat{h}(\xi_2,\tau_2)
    \nonumber\\
    &\lesssim&
    \|f\|_{X_{0,\frac{1}{2}}}\>
    \|g\|_{X_{0,\frac{1}{2}}}\>
    \|h\|_{X_{0,\frac{1}{2}}},
\end{eqnarray*}
for any $f, g, h\in X_{0,\frac{1}{2}}$ such that supp$_\xi
\hat{g}\in [0,N]$ and supp$_\xi \hat{h}\in [N,+\infty)$ belongs to
$[N,+\infty)$, thus $|\xi_1|\leq|\xi_2|$. Note that $|\xi|\sim
|\xi_2|$ and $\langle\xi\rangle^{s}-\langle\xi_2\rangle^{s} \lesssim
|\xi_1||\xi_2|^{s-1}$, the left-hand side is reduced to
\begin{equation}\label{2.13}
\displaystyle\int_{\ast}
|\xi|^{\frac{1}{2}}|\xi_1|^{\frac{1}{2}}|\xi_2|^{\frac{1}{2}}
\hat{f}(\xi,\tau)\>\hat{g}(\xi_1,\tau_1)\>\hat{h}(\xi_2,\tau_2).
\end{equation}
By the fact that
$$
\tau-\xi^3-(\tau_1-\xi_1^3)-(\tau_2-\xi_2^3)=-3\xi\xi_1\xi_2,
$$
we may assume that $|\tau-\xi^3|\gtrsim |\xi\xi_1\xi_2|$ by
symmetry. Therefore we reduce to control
\begin{equation}\label{5.7}
\displaystyle\int_{\ast} \langle\tau-\xi^3\rangle^{\frac{1}{2}}
\hat{f}(\xi,\tau)\>\hat{g}(\xi_1,\tau_1)\>\hat{h}(\xi_2,\tau_2).
\end{equation}
By H\"{o}lder's inequality and Lemma 2.1, (\ref{5.7}) is
bounded by
\begin{equation}\label{215}
\|f\|_{X_{0,\frac{1}{2}}}\> \|g\|_{X_{0,\frac{1}{3}}}\>
\|h\|_{X_{0,\frac{1}{3}}}.
\end{equation}
This completes the proof of the lemma. \hfill$\Box$
\end{proof}

{\it \noindent Remark.} By (\ref{5.8}) and the proof above (see (\ref{215})), we
actually have
\begin{equation}
\left\|\varphi\>\partial_x[J_x^s,u]v\right\|_{X_{\frac{1}{2},-\frac{1}{2}}}
    \lesssim
    \|u\|_{X_{\frac{1}{2},\frac{1}{2}}}\>
    \|v\|_{X_{s-\frac{1}{2},\frac{1}{2}}},\label{L5.4'}
\end{equation}
for any $\varphi\in H^{\frac{1}{2}-}(\R)$ and $u,v$
as Lemma 2.4.

We note that Lemma 2.4 gives a good estimate under the restrictions to
$u,v$ on the space-frequency space. A nature question is: what happens if
we drop the restrictions? Now we give a general result about it.
\begin{lem}
For any $s\in \R$, if the functions $u,v$ are zero
$x-$mean for all $t$, then
\begin{equation}
\left\|\partial_x[J_x^s,u]v\right\|_{X_{s_1,-\frac{1}{2}}}
    \lesssim
    \|u\|_{X_{s_2,\frac{1}{2}}}\>
    \|v\|_{X_{s_3,\frac{1}{2}}},\label{L2.5}
\end{equation}
where $s_1,s_2,s_3$ satisfy
\begin{equation}\label{c}
s_1=-\frac{1}{2}+a-c;\quad
s_2=-\frac{1}{2}+s+a+b-c;\quad
s_3=-\frac{1}{2}+c-b,
\end{equation}
for any $a,b,c\geq 0$ and $a+b-2c\leq 1-s$.
\end{lem}
\begin{proof}
By duality and Plancherel's identity, it suffices to show
\begin{eqnarray*}
    &&\displaystyle\int_{\ast}
    \xi\langle\xi\rangle^{s_1}
    \langle\xi_1\rangle^{-s_2}
    \langle\xi_2\rangle^{-s_3}
    (\langle\xi\rangle^{s}-\langle\xi_2\rangle^{s})\>
    \hat{f}(\xi,\tau)\>\hat{g}(\xi_1,\tau_1)\>\hat{h}(\xi_2,\tau_2)
    \nonumber\\
    &\lesssim&
    \|f\|_{X_{0,\frac{1}{2}}}\>
    \|g\|_{X_{0,\frac{1}{2}}}\>
    \|h\|_{X_{0,\frac{1}{2}}},
\end{eqnarray*}
for any $f, g, h\in
X_{0,\frac{1}{2}}$.

First, if $|\xi|\sim |\xi_2| $, then $|\xi_1|\lesssim |\xi|$. The
same proof as Lemma 2.4, the left-hand side is reduced to
$$
\displaystyle\int_{\ast}
    |\xi|^{s_1-s_3+s}
    |\xi_1|^{1-s_2}
    \>
    \hat{f}(\xi,\tau)\>\hat{g}(\xi_1,\tau_1)\>\hat{h}(\xi_2,\tau_2).
$$
It encounters (\ref{2.13}) when
\begin{equation}\label{c1}
s_1-s_3+s-1\leq 0, \mbox{\quad and\quad}
s_1-s_2-s_3+s-\frac{1}{2}\leq 0.
\end{equation}
Thus we have the claimed result.

Second, if $|\xi|\ll |\xi_2|\sim |\xi_1| $, then the left-hand side is reduced to
$$
\displaystyle\int_{\ast}
    |\xi|^{1+s_1}
    |\xi_1|^{-s_2-s_3+s}
    \>
    \hat{f}(\xi,\tau)\>\hat{g}(\xi_1,\tau_1)\>\hat{h}(\xi_2,\tau_2).
$$
It encounters (\ref{2.13}) again when
\begin{equation}\label{c2}
-s_2-s_3+s-1\leq 0, \mbox{\quad and\quad} s_1-s_2-s_3+s-\frac{1}{2}\leq 0.
\end{equation}

Third, if $|\xi|\sim |\xi_1|\gg |\xi_2|$, then the left-hand side is reduced to
$$
\displaystyle\int_{\ast}
    |\xi|^{1+s_1-s_2+s}
    |\xi_2|^{-s_3}
    \>
    \hat{f}(\xi,\tau)\>\hat{g}(\xi_1,\tau_1)\>\hat{h}(\xi_2,\tau_2).
$$
It also meets (\ref{2.13}) when
\begin{equation}\label{c3}
s_1-s_2+s\leq 0, \mbox{\quad and\quad} s_1-s_2-s_3+s-\frac{1}{2}\leq 0.
\end{equation}

We set that
$$
s_2+s_3-s+1=a;\quad -s_1+s_2-s=b;\quad -s_1+s_2+s_3-s+\frac{1}{2}=c,
$$
then inserting them into (\ref{c1})--(\ref{c3}), we give the
condition (\ref{c}) in the lemma. \hfill$\Box$
\end{proof}

As a consequence of Lemma 2.4, we have
\begin{cor}
For any $s\leq 0$, if the functions $u,v$ are zero
$x-$mean for all $t$, then
\begin{equation}
\left\|\varphi\>\partial_x[J_x^s,u]v\right\|_{X_{\frac{1}{2}+s,-\frac{1}{2}}}
    \lesssim
    \|u\|_{X_{\frac{1}{2}+2s,\frac{1}{2}}}\>
    \|v\|_{X_{-\frac{1}{2},\frac{1}{2}}},\label{L5.4''}
\end{equation}
for any $\varphi\in H^{\frac{1}{2}-}(\R)$.
\end{cor}
\begin{proof} We choose $a=1+s,b=c=0$ in (\ref{L5.4'}), and use
(\ref{5.8}) to get the result.
\hfill$\Box$
\end{proof}

At the end of this section, we give a special bilinear estimate as
follows, which might be not sharp but enough to use in the
article. A similar one has been presented in \cite{Ts2}.
\begin{lem}
For any $s>-\dfrac{1}{4}$, $u\in X_{s,\frac{1}{2}}^\delta, a\in
H^s(\T)$ such that $u,a$ are zero $x-$mean, then
\begin{equation}\label{L25}
 \|\psi(t/\delta)\>\partial_x(ua)\|_{Z^{s}}
 \lesssim
 \delta^{\frac{1}{3}-}\|u\|_{X_{s,\frac{1}{2}}^\delta}
 \|a\|_{H^s}.
\end{equation}
\end{lem}
\begin{proof}
According to the two contributions to the $Z^s$-norm, we need to
show
\begin{eqnarray}
\|\psi(t/\delta)\>\partial_x(ua)\|_{X_{s,\frac{1}{2}}} \lesssim
 \delta^{\frac{1}{3}-}\|u\|_{X_{s,\frac{1}{2}}^\delta}
 \|a\|_{H^s},\label{2.15}\\
\left\|\dfrac{\xi[\psi(\cdot/\delta)\ast\hat{u}\ast \hat{a}]}
{\langle\tau-\xi^3\rangle}\right\|_{L^2_\xi L^1_\tau} \lesssim
 \delta^{\frac{1}{3}-}\|u\|_{X_{s,\frac{1}{2}}^\delta}\|a\|_{H^s}.
 \label{2.16}
\end{eqnarray}

For (\ref{2.15}), we only need to show
\begin{equation}\label{217}
\displaystyle\int_{\ast} \frac{\xi\langle\xi\rangle^s}
{\langle\xi_1\rangle^s\langle\xi_2\rangle^s}
\hat{f}(\xi,\tau)\hat{g}(\xi_1,\tau_1)\widehat{h_\delta}(\xi_2,\tau_2)
\lesssim \delta^{\frac{1}{3}-}\|f\|_{X_{0,\frac{1}{2}}^\delta}
\>\|g\|_{X_{0,\frac{1}{2}}^\delta}\>\|h\|_{L^2_{x}}
\end{equation}
for any $f,g\in X_{0,\frac{1}{2}}^\delta$, $h\in L^2_{xt}$, where we
denote $h_\delta(x,t)=\psi(t/\delta)h(x)$. Moreover, we may assume
that $f=\psi(t/\delta)f,g=\psi(t/\delta)g$.

Note that
$$
|(\tau-\xi^3)-(\tau_1-\xi_1^3)-\tau_2|
=|\xi_2(\xi^2+\xi\xi_1+\xi_1^2)| \geq |\xi\xi_1\xi_2|,
$$
we may split the integral (\ref{217}) into three parts:
$$
{\bf (a)}\,|\tau-\xi^3|\gtrsim |\xi||\xi_1||\xi_2|;\ \
{\bf (b)}\,|\tau_1-\xi_1^3|\gtrsim |\xi||\xi_1||\xi_2|;\ \
{\bf (c)}\,|\tau_2|\gtrsim |\xi||\xi_1||\xi_2|.
$$

${\bf (a)}:|\tau-\xi^3|\gtrsim |\xi||\xi_1||\xi_2|$, we split the
integral again into two parts: {\bf Part 1:} $|\xi|\lesssim
|\xi_1|$; {\bf Part 2:} $|\xi_1|\ll |\xi| \sim |\xi_2|$.

{\bf Part 1:} $|\xi|\lesssim |\xi_1|$. The left-hand side of
(\ref{217}) restricted in this part is controlled by
\begin{eqnarray*}
&&
\displaystyle\int_{\ast}|\xi_2|^{-\frac{1}{2}-s}
\langle\tau-\xi^3\rangle^{\frac{1}{2}}
\hat{f}(\xi,\tau)\hat{g}(\xi_1,\tau_1)\widehat{h_\delta}(\xi_2,\tau_2)\\
&\lesssim& \|f\|_{X_{0,\frac{1}{2}}}\>\|g\|_{L^4_{xt}}
\>\left\|J_x^{-\frac{1}{2}-s}h_\delta\right\|_{L^4_{xt}}.
\end{eqnarray*}
Since $s\geq -\dfrac{1}{4}$, by Sobolev's inequality, we have,
\begin{eqnarray*}
\left\|J_x^{-\frac{1}{2}-s}h_\delta\right\|_{L^4_{xt}}
=\|\psi(t/\delta)\|_{L^4_t}\>\left\|J_x^{-\frac{1}{2}-s}h\right\|_{L^4_x}
\lesssim \delta^{\frac{1}{4}}\>\|h\|_{L^2_x},
\end{eqnarray*}
where we note that for any $1\leq p\leq \infty$,
\begin{equation}\label{phi}
\|\psi(\cdot/\delta)\|_{L^p}= \delta^{\frac{1}{p}}\|\psi\|_{L^p}.
\end{equation}
Then (\ref{217}) restricted in this part follows by (\ref{L5.1}) and
(\ref{delta}).

{\bf Part 2:} $|\xi_1|\ll |\xi| \sim |\xi_2|$. The left-hand side of
(\ref{217}) restricted in this part is bounded by
\begin{eqnarray*}
&&
\displaystyle\int_{\ast}|\xi_1|^{-\frac{1}{2}-s}
\langle\tau-\xi^3\rangle^{\frac{1}{2}}
\hat{f}(\xi,\tau)\hat{g}(\xi_1,\tau_1)\widehat{h_\delta}(\xi_2,\tau_2)\\
&\lesssim& \|f\|_{X_{0,\frac{1}{2}}}
\>\left\|J_x^{-\frac{1}{2}-s}g\right\|_{L^4_{t}L^\infty_{x}}
\>\|h_\delta\|_{L^4_{t}L^2_{x}}.
\end{eqnarray*}
Further, since $s>- \dfrac{1}{4}$, by Sobolev's inequality, $
\left\|J_x^{-\frac{1}{2}-s}g\right\|_{L^4_{t}L^\infty_{x}}
\lesssim\|g\|_{L^4_{xt}}, $ and $
\|h_\delta\|_{L^4_{t}L^2_{x}}\lesssim
\delta^{\frac{1}{4}}\>\|h\|_{L^2}. $ Then (\ref{217}) restricted in
this part follows by (\ref{L5.1}) and (\ref{delta}) again.

${\bf (b)}:|\tau_1-\xi_1^3|\gtrsim |\xi||\xi_1||\xi_2|$, we also
split the integral again into two parts: {\bf Part 1:}
$|\xi|\lesssim |\xi_1|$; {\bf Part 2:} $|\xi_1|\ll |\xi| \sim
|\xi_2|$. But \textbf{Part 1} is similar to Part (a)(1), so we only
consider Part 2.

{\bf Part 2:} $|\xi_1|\ll |\xi| \sim |\xi_2|$. The left-hand side of
(\ref{217}) restricted in this part is controlled by
\begin{eqnarray*}
&&
\displaystyle\int_{\ast}|\xi_1|^{-\frac{1}{2}-s}
\langle\tau_1-\xi_1^3\rangle^{\frac{1}{2}}
\hat{f}(\xi,\tau)\hat{g}(\xi_1,\tau_1)\widehat{h_\delta}(\xi_2,\tau_2)\\
&\lesssim& \|f\|_{L^4_{xt}}
\>\left\|J_x^{-\frac{1}{2}-s}\Gamma g\right\|_{L^2_{t}L^4_{x}}
\>\|h_\delta\|_{L^4_{t}L^2_{x}},
\end{eqnarray*}
where the operator $\Gamma$ defined by $\widehat{\Gamma g}(\xi,\tau)
=\langle\tau-\xi^3\rangle^{\frac{1}{2}}\hat{g}(\xi,\tau)$. By
Sobolev's inequality, $ \left\|J_x^{-\frac{1}{2}-s}\Gamma
g\right\|_{L^2_{t}L^4_{x}} \lesssim\|\Gamma
g\|_{L^2_{xt}}=\|g\|_{X_{0,\frac{1}{2}}} $ since $s\geq
-\dfrac{1}{4}$. Then (\ref{217}) follows by (\ref{L5.1}),
(\ref{phi}) and (\ref{delta}).

${\bf (c)}:|\tau_2|\gtrsim |\xi||\xi_1||\xi_2|$, then the left-hand
side of (\ref{217}) is bounded by
\begin{eqnarray*}
\displaystyle\int_{\ast} |\tau_2|^{\frac{1}{2}}
\hat{f}(\xi,\tau)\hat{g}(\xi_1,\tau_1)\widehat{h_\delta}(\xi_2,\tau_2)
\lesssim \|f\|_{L^4_{xt}}\>\|g\|_{L^4_{xt}}
\>\left\|D_t^{\frac{1}{2}}h_\delta\right\|_{L^2_{xt}} \lesssim
\delta^{\frac{1}{3}-}\|f\|_{X_{0,\frac{1}{2}}}
\>\|g\|_{X_{0,\frac{1}{2}}}\>\|h\|_{L^2_{x}}
\end{eqnarray*}
by (\ref{L5.1}) and the fact
$\big\|D^{\frac{1}{2}}\psi(\cdot/\delta)\big\|_{L^2}=
\big\|D^{\frac{1}{2}}\psi\big\|_{L^2}$. Thus we prove the result
(\ref{217}).

For (\ref{2.16}), if $|\tau_1-\xi_1^3|\gtrsim |\xi||\xi_1||\xi_2|$
or $|\tau_2|\gtrsim |\xi||\xi_1||\xi_2|$, then by H\"{o}lder's
inequality, (\ref{2.16}) is sufficient if
\begin{equation}\label{2.26}
\left\|\dfrac{\xi[\psi(\cdot/\delta)\ast\hat{u}\ast \hat{a}]}
{\langle\tau-\xi^3\rangle^{\frac{1}{2}-}}\right\|_{L^2_{\xi\tau}}
\lesssim
 \delta^{\frac{1}{3}-}\|u\|_{X_{s,\frac{1}{2}}^\delta}\|a\|_{H^s}.
\end{equation}
But it can be shown as (\ref{2.15}) above. Therefore, we only
consider
$$
|\tau-\xi^3|\gtrsim |\xi||\xi_1||\xi_2|.
$$

Moreover, in the event that
$$
|\tau_1-\xi_1^3|\gtrsim |\xi\xi_1\xi_2|^{0+},
$$
(\ref{2.16}) is reduced to (\ref{2.26}) and turned further to
\begin{eqnarray*}
&&\displaystyle\int_{\ast} |\xi|^{\frac{1}{2}+s}|\xi_1|^{-\frac{1}{2}-s}
|\xi_2|^{-\frac{1}{2}-s}\langle\tau_1-\xi_1^3\rangle^{0+}
\hat{f}(\xi,\tau)\hat{g}(\xi_1,\tau_1)\widehat{h_\delta}(\xi_2,\tau_2)\\
&\lesssim&
 \delta^{\frac{1}{3}-}\|f\|_{L^2_{xt}}
 \|g\|_{X_{0,\frac{1}{2}}}\|h\|_{L^2_{x}},
\end{eqnarray*}
which can be shown by the same argument used in Part (a) above.
The event that $|\tau_2|\gtrsim |\xi\xi_1\xi_2|^{0+}$ is similar.
Therefore, we only need to consider the case that
$$
\tau_1-\xi_1^3=-\xi_2(\xi^2+\xi\xi_1+\xi_1^2)+O(|\xi\xi_1\xi_2|^{0+}).
$$
Let the set
$$
\Omega(\xi)=\{\eta\in\R: \eta=-\xi_2(\xi^2+\xi\xi_1+\xi_1^2)
+O(|\xi\xi_1\xi_2|^{0+}) \mbox{ for any } \xi_1, \xi_2\in \Z \mbox{ with }
\xi=\xi_1+\xi_2\},
$$
then similar to Lemma 7.6 in \cite{CKSTT2}, we have
$$
|\Omega(\xi)\cap \{\eta: |\eta|\sim M\}|\lesssim M^{\frac{2}{3}},
$$
which leads to
\begin{equation}\label{219}
\left(\displaystyle\int \langle\tau-\xi^3\rangle^{-1}
\chi_{\Omega(\xi)}(\tau-\xi^3)\, d\tau\right)^{\frac{1}{2}}
\lesssim
1.
\end{equation}
Indeed, it can be easily proved
by the dyadic decomposition to the integration.

Using (\ref{219}) and H\"{o}lder's inequality, the left-hand side of
(\ref{2.16}) is controlled by
$$
\left\|\dfrac{\xi[\psi(\cdot/\delta)\ast\hat{u}\ast \hat{a}]}
{\langle\tau-\xi^3\rangle^\frac{1}{2}}\right\|_{L^2_{\xi\tau}} .
$$
Again  (\ref{2.16}) follows from the same argument used in Part
(a).
 \hfill$\Box$
\end{proof}

{\it \noindent Remark.} By (\ref{5.8}) and the proof above, we
actually have
\begin{equation*}
 \|\varphi\psi(t/\delta)\>\partial_x(ua)\|_{Z^{s}}
 \lesssim
 \delta^{\frac{1}{3}-}\|u\|_{X_{s,\frac{1}{2}}^\delta}
 \|a\|_{H^s},
\end{equation*}
for any $\varphi\in H^{\frac{1}{2}-}(\R)$ and $u,v$ as Lemma 2.7.
As a consequence, we have
\begin{equation}\label{L25'}
 \|\varphi\psi(t/\delta)\>\partial_x(ua)\|_{X_{s+l,-\frac{1}{2}}}
 \lesssim
 \|u\|_{X_{s+l,\frac{1}{2}}^\delta}
 \|a\|_{H^s}+
 \|u\|_{X_{s,\frac{1}{2}}^\delta}
 \|a\|_{H^{s+l}},
\end{equation}
for any $s>-\dfrac{1}{4},l\geq 0.$

 \vspace{0.3cm}
\section{Well-posedness and Bounded Absorbing Sets}

In this section, we first give the local well-posedness,
then we apply the I-method and the multilinear correction technique
to prove the global well-posedness,
finally we obtain bounded absorbing sets.

\subsection{The Local Well-posedness}
Compared to the KdV equation, the equation (\ref{KdV}) lacks the
solution of scale invariance, so we have to dig some additional
factors of $\delta$ from the estimates, which is of importance for
us. The first one is a refined local result which is a contrast to
Proposition 4 in \cite{CKSTT2}.

By employing the bilinear estimates (\ref{L21}) (which replace
(7.33) in \cite{CKSTT2}), we give the following local result (which
instead of Proposition 4 in \cite{CKSTT2}).
\begin{prop}
Let $s\geq -\dfrac{1}{2}$, $I=I_{N,s}$, $f\in \dot{H}^s(\T)$, then
(\ref{KdV}) (\ref{initial data}) are locally well-posed for the
initial data $u_0$ satisfies $Iu_0\in \dot{L}^2(\T)$, with the lifetime
$\delta$ satisfying
\begin{equation}
 \delta\sim (\|I_{N,s}u_0\|_{L^2}+\|I_{N,s}f\|_{L^2})^{-3-}.
 \label{Edelta}
\end{equation}
Further, the solution satisfies the estimate
\begin{equation}
    \|I_{N,s}u\|_{Y^0}
    \lesssim
    \|I_{N,s}u_0\|_{L^2}+\|I_{N,s}f\|_{L^2}.\label{LSE}
\end{equation}
\end{prop}

{\it \noindent  Remark.} The improvement in this local result is to
give a refined estimates on lifetime $\delta$, which strongly effect
on the global well-posedness of the weakly damped, forced KdV
equation who lacks the solution of  scale invariance (see the
following subsection for more details).

By (\ref{LSE}), we have further the control of the solution as
\begin{equation}\label{LSE2}
\sup\limits_{t\in [t_0-\delta,t_0+\delta]} \|I_{N,s}u(t)\|_{L^2_x}
\lesssim \|I_{N,s}u(t_0)\|_{L^2}+\|I_{N,s}f\|_{L^2},
\end{equation}
if we take ``$t_0$'' for the initial time.
 \vspace{0.3cm}

\subsection{The Global Well-posedness}
Now we are further to consider the global well-posedness and the
existence of bounded absorbing sets.
The argument here is mainly
the I-method in \cite{CKSTT2}. However, some estimates
and the iteration process used in \cite{CKSTT2}
should be rebuilt. We show that, due to the presence of the weak damping,
the energy will not increase unlimitedly.

First we define the symmetrization of a $k-$multipler
$m:\Z^k\rightarrow \R$ by
$$
[m]_{sym(\xi)}=\frac{1}{n!}\sum_{\sigma\in S_n}m(\sigma(\xi)),
$$
where $S_n$ is the group of all permutations on $n$ objects. We say
$m$ is symmetric if $m(\xi)=m(\sigma(\xi))$.

Define the $k-$multiplier
$$
  \Lambda_k(m; u_1,\cdots,u_k)
 =\displaystyle\int_{\sum\limits_{j=1}^k\xi_j=0}m(\xi_1,\cdots,\xi_k)
   \prod_{j=1}^k\F_xu_j(\xi_j,t)\,(d\xi_1)_\T\cdots (d\xi_{k-1})_\T.
$$
We write $\Lambda_k(m)=\Lambda_k(m; u,\cdots,u)$ for short. Then by
a direct computation, we have
\begin{equation}\label{2.9}
\begin{array}{c}
\dfrac{d}{dt}\Lambda_k(m)=-k\gamma \Lambda_k(m)+\Lambda_k(\alpha_k
m)
+k\Lambda_k(m;u,\cdots,u,f)\\
-\dfrac{i}{2}k\Lambda_{k+1}(m(\xi_1,\cdots,\xi_{k-1},\xi_{k}+\xi_{k+1})
(\xi_{k}+\xi_{k+1})),
\end{array}
\end{equation}
where the multiplier $m$ is symmetric and
$$
\alpha_k\equiv i(\xi_1^3+\cdots+\xi_k^3).
$$
Note that the fourth term of (\ref{2.9}) may be symmetrized.

Now let $m(\xi_1,\xi_2)=m(\xi_1)m(\xi_2)$, denote the modified
energy as
$$
E^2_I(t)\equiv\|Iu(t)\|_{L^2}^2=\Lambda_2(m),
$$
then by (\ref{2.9}) and note that $\alpha_2(\xi_1,\xi_2)=0$ when
$\xi_1+\xi_2=0$, we have
$$
\dfrac{d}{dt}E^2_I(t) =-2\gamma
E^2_I(t)+2\Lambda_2(m;u,f)+\Lambda_3(M_3),
$$
where
$$
M_3(\xi_1,\xi_2,\xi_3)
 =-\frac{i}{3}\big(m^2(\xi_1)\xi_1+m^2(\xi_2)\xi_2+m^2(\xi_3)\xi_3\big).
$$
Define a new modified energy $E^3_I(t)$ by \beq
E^3_I(t)=\Lambda_3(\sigma_3)+E^2_I(t), \label{E3} \eeq where
$$
\sigma_3=-\dfrac{M_3}{\alpha_3}.
$$
Then one has
$$
\dfrac{d}{dt}E^3_I(t)=-\gamma E^3_I(t)-\gamma E^2_I(t)-2\gamma\Lambda_3(\sigma_3)+2\Lambda_2(m;u,f)+
3\Lambda_3(\sigma_3;u,u,f)+ \Lambda_4(M_4),
$$
where
$$
M_4=-\dfrac{3}{2}i[\sigma_3(\xi_1,\xi_1,\xi_3+\xi_4)(\xi_3+\xi_4)]_{sym}.
$$

Define another new modified energy $E^4_I(t)$ again by \beq
E^4_I(t)=\Lambda_4(\sigma_4)+E^3_I(t), \label{E4} \eeq with
$$
\sigma_4=-\dfrac{M_4}{\alpha_4}.
$$
Then one has
$$
\begin{array}{c}
\dfrac{d}{dt}E^4_I(t)=-\gamma E^4_I(t)-\gamma E^2_I(t)
-2\gamma\Lambda_3(\sigma_3)-3\gamma\Lambda_4(\sigma_4)
+2\Lambda_2(m;u,f)
\\
+3\Lambda_3(\sigma_3;u,u,f)+4\Lambda_4(\sigma_4;u,u,u,f)+
\Lambda_5(M_5),
\end{array}
$$
where
$$
M_5=-2i[\sigma_4(\xi_1,\xi_1,\xi_3,\xi_4+\xi_5)(\xi_4+\xi_5)]_{sym}.
$$
Therefore,
\begin{equation}\label{2.12}
    E^4_I(t+\delta)=e^{-\gamma\delta}E^4_I(t)
    +\displaystyle\int_0^\delta\!
    e^{-\gamma(\delta-t')}[F(t'+t)-\gamma E^2_I(t'+t)]\,dt',
\end{equation}
where
$$
\begin{array}{c}
F=-2\gamma\Lambda_3(\sigma_3)-3\gamma\Lambda_4(\sigma_4)
+2\Lambda_2(m;u,f)\\
\qquad\qquad\qquad+3\Lambda_3(\sigma_3;u,u,f)+4\Lambda_4(\sigma_4;u,u,u,f)+
\Lambda_5(M_5).
\end{array}
$$

First, by (\ref{LSE2}), we have
\begin{equation}\label{2.22}
\gamma \displaystyle\int_0^\delta\!
    e^{-\gamma(\delta-t')}E^2_I(t'+t)\,dt'
\geq c\gamma\delta \|Iu(t)\|_{L^2}^2-C\gamma\delta\|If\|_{L^2}^2,
\end{equation}
for some small $c$, large $C>0$.

Now we give some multilinear estimates on $\Lambda_k$. The first
one is an improvement of the results in Lemma 6.1 in
\cite{CKSTT2}, which is of importance in this situation.
\begin{lem}
Let $s\geq-\dfrac{3}{4},I=I_{N,s}$,  then
\begin{eqnarray}
|\Lambda_3(\sigma_3;u_1,u_2,u_3)| &\lesssim&
N^{-\frac{3}{2}}\prod\limits_{j=1}^3\|Iu_j(t)\|_{L^2};\label{A3}\\
|\Lambda_4(\sigma_4;u_1,u_2,u_3,u_4)| &\lesssim&
N^{-3}\prod\limits_{j=1}^4\|Iu_j(t)\|_{L^2}.\label{A4}
\end{eqnarray}
\end{lem}
\begin{proof}
For (\ref{A3}), since $\xi_1+\xi_2+\xi_3=0$, by symmetry we may
assume $|\xi_1|\sim|\xi_2|\geq |\xi_3|$. Note that $\sigma_3$
vanishes when $|\xi_j|\leq N$ for $j=1,2,3$, so we may assume
further that $|\xi_1|,|\xi_2|\gtrsim N$.

Set
$$
\triangle\equiv\dfrac{|\sigma_3|}{m(\xi_1)m(\xi_2)m(\xi_3)}
=\dfrac{|M_3(\xi_1,\xi_2,\xi_3)|}
       {|\alpha_3(\xi_1,\xi_2,\xi_3)|\,m(\xi_1)m(\xi_2)m(\xi_3)},
$$
then (\ref{A3}) follows if we show
$$
|\Lambda_3(\triangle;u_1,u_2,u_3)|\lesssim
N^{-\frac{3}{2}}\|u_1\|_{L^2}\|u_2\|_{L^2} \|u_3\|_{L^2}.
$$
Note that
$$
|M_3(\xi_1,\xi_2,\xi_3)|\lesssim m^2(\xi_{min})|\xi_{min}|, \quad
\alpha_3(\xi_1,\xi_2,\xi_3)=3\xi_1\xi_2\xi_3,
$$
thus, noting that $s\geq -\dfrac{3}{4}$, we have
$$
\triangle\lesssim\dfrac{1}{|\xi_1\xi_2|\,m(\xi_1)m(\xi_2)}\sim
N^{2s}|\xi_1|^{-1-s}|\xi_2|^{-1-s}\lesssim
N^{-\frac{3}{2}}|\xi_1|^{-\frac{1}{4}}|\xi_2|^{-\frac{1}{4}}.
$$
Therefore, by H\"{o}lder and Sobolev's inequalities,  we have \beqs
\begin{split}
 |\Lambda_3(\triangle;u_1,u_2,u_3)|&\lesssim
N^{-\frac{3}{2}}|\Lambda_3(|\xi_1|^{-\frac{1}{4}}|\xi_2|^{-\frac{1}{4}};u_1,u_2,u_3)|\\
&\lesssim
N^{-\frac{3}{2}}\left\|D_x^{-\frac{1}{4}}u_1\right\|_{L^4}\left\|D_x^{-\frac{1}{4}}u_2\right\|_{L^4}
\|u_3\|_{L^2}\\
&\lesssim N^{-\frac{3}{2}}\|u_1\|_{L^2}\|u_2\|_{L^2} \|u_3\|_{L^2}.
\end{split}
\eeqs

Now we turn to (\ref{A4}), and set
$$
\tilde{\triangle}\equiv\dfrac{|\sigma_4|}{\prod\limits_{j=1}^4m(\xi_j)},
$$
then (\ref{A4}) suffices if we show
$$
|\Lambda_4(\tilde{\triangle};u_1,u_2,u_3,u_4)| \lesssim
N^{-3}\prod\limits_{j=1}^4\|u_j(t)\|_{L^2}.
$$

By Lemma 4.4 in \cite{CKSTT2}, we have
$$
\tilde{\triangle}\lesssim
\dfrac{1}{\prod\limits_{j=1}^4m(\xi_j)\max\{N,|\xi_j|\}}.
$$
Since $m(\xi)\sim N^{-s}\max\{N,|\xi|\}^s$, and noting that $s\geq
-\dfrac{3}{4}$, we have
$$
\tilde{\triangle}\lesssim
N^{-3}\prod\limits_{j=1}^4|\xi_j|^{-\frac{1}{4}}.
$$
Therefore, we control
$|\Lambda_4(\tilde{\triangle};u_1,u_2,u_3,u_4)|$ by
$$
N^{-3}\prod\limits_{j=1}^4\left\|D^{-\frac{1}{4}}u_j(t)\right\|_{L^4_x}
\lesssim N^{-3}\prod\limits_{j=1}^4\|u_j(t)\|_{L^2},
$$
by Sobolev's inequality. \hfill$\Box$
\end{proof}

\begin{lem}
Let $s\geq-\dfrac{1}{2},I=I_{N,s}, \beta=\dfrac{5}{2}-,
0<\delta\leq 1$, then
\begin{eqnarray}
 \Big|\displaystyle\int_0^\delta
 \Lambda_5(M_5)\,dt\Big|
 \lesssim
 N^{-\beta}\delta^{\frac{1}{3}-}\|Iu\|_{Y^0}^5.
 \label{A5}
\end{eqnarray}
\end{lem}
\begin{proof}
We only give the modification of the proof in \cite{CKSTT2} here.
By the modified bilinear estimate (\ref{L21}), we replace (8.1) in
\cite{CKSTT2} by
$$
\displaystyle\left|\int_0^\delta\!\!\int_\T \prod\limits_{j=0}^5
w_j(x,t)\,dx dt\right| \lesssim
\delta^{\frac{1}{3}-}\prod\limits_{j=0}^3
\|w_j\|_{Y^{\frac{1}{2}}}\>\|w_4\|_{X_{-\frac{1}{2},\frac{1}{2}}}\>
\|w_5\|_{X_{-\frac{1}{2},\frac{1}{2}}},
$$
where $w_j(x,t)$ are $\T$-periodic, and zero $x-$mean for any $t$.
Based on that, we lead to the result claimed. \hfill$\Box$
\end{proof}

Now we consider {\it a priori} estimates of the solutions, and from
now on we set $\delta\sim N^{3s-}$ by (\ref{Edelta}) (the large
number $N$ will be chosen later). Back to (\ref{2.12}), we apply
(\ref{A3})--(\ref{A5}) to yield a bound of
$\displaystyle\int_0^\delta\! F(t'+t)\,dt'$ as
\begin{eqnarray}
 &\delta\sup\limits_{t'\in [t-\delta,t+\delta]}
 \Big(N^{-\frac{3}{2}}\>\|Iu(t')\|_{ L^2_x}^3
 + N^{-3}\>\|Iu(t')\|_{L^2_x}^4
 + \>\|Iu(t')\|_{L^2_x}\>\|If\|_{L^2_x}\nonumber\\
 &+ N^{-\frac{3}{2}}\>\|Iu(t')\|_{L^2_x}^2\>\|If\|_{L^2_x}
 + N^{-3}\>\|Iu(t')\|_{L^2_x}^3\>\|If\|_{L^2_x}\Big)
 + N^{-\beta}\delta^{\frac{1}{3}-}
 \>\|Iu\|_{Y^0}^5.\label{2.17}
\end{eqnarray}
By (\ref{II}) and recalling that $\delta\sim N^{3s-}$, we have
\begin{equation}\label{2.20}
\delta N^{-\frac{3}{2}}\|If\|_{L^2_x}^3,\,\,
\delta N^{-3}\|If\|_{L^2_x}^4,\,\,
\delta\|If\|_{L^2_x}^2,\,\,
N^{-\beta}\delta^{\frac{1}{3}-}\|If\|_{L^2_x}^5
\leq K(f),
\end{equation}
for some positive constant $K$ only dependent on $\|f\|_{H^s}$.
Therefore, by (\ref{LSE}), (\ref{LSE2}) and (\ref{2.20}),
(\ref{2.17}) is controlled by
\begin{eqnarray}
 &K(f)
 +\dfrac{c}{2}\gamma\delta\>\|Iu(t)\|_{L^2_x}^2
 +C\delta\|Iu(t)\|_{L^2_x}^2 \Big(N^{-\frac{3}{2}}\>\|Iu(t)\|_{L^2_x}
 +N^{-3}\>\|Iu(t)\|_{L^2_x}^2
  \nonumber\\
 &
 +N^{-\frac{3}{2}}\>\|If\|_{L^2_x}
 +N^{-3}\>\|Iu(t)\|_{L^2_x}\>\|If\|_{L^2_x}\Big)
 +CN^{-\beta}\delta^{\frac{1}{3}-}\|Iu(t)\|_{L^2_x}^5
 \label{2.19}
\end{eqnarray}
for some constant $C$ independent on $N,u$ and $f$.  Then by
(\ref{II}) again, we control (\ref{2.19}) and finally
get
\begin{eqnarray}
 \displaystyle\int_0^\delta\! F(t'+t)\,dt'
 &\lesssim&
 K(f)+\dfrac{c}{2}\gamma\delta\|Iu(t)\|_{L^2_x}^2
 +C\delta\|Iu(t)\|_{L^2_x}^2\cdot
 \nonumber\\
 &&
 \cdot\big(N^{-\frac{3}{2}}\>\|Iu(t)\|_{L^2_x}
 +N^{-3}\>\|Iu(t)\|_{L^2_x}^2+N^{-\frac{3}{2}-s}\>\|f\|_{H^s}
 \nonumber\\
 &&
 +N^{-3-s}\>\|Iu(t)\|_{L^2_x}\>\|f\|_{H^s}
 +N^{-\beta}\delta^{-\frac{2}{3}+}\>\|Iu(t)\|_{L^2_x}^3\big).\label{2.21}
\end{eqnarray}

Combining with (\ref{2.12}), (\ref{2.22}), (\ref{2.20}), (\ref{2.21}),
and choosing $\epsilon$ small enough, we have
\begin{eqnarray}
    E^4_I(t+\delta)
&\leq&
    e^{-\gamma\delta}E^4_I(t)
    +K(f)+\delta\|Iu(t)\|_{L^2_x}^2
    \Big(CN^{-\frac{3}{2}}\>\|Iu(t)\|_{L^2_x}
    +CN^{-3}\>\|Iu(t)\|_{L^2_x}^2\nonumber\\
&&
    +CN^{-\frac{3}{2}-s}\>\|f\|_{H^s}
    +CN^{-3-s}\>\|Iu(t)\|_{L^2_x}\>\|f\|_{H^s}\nonumber\\
&&
    +CN^{-\beta}\delta^{-\frac{2}{3}+}\>\|Iu(t)\|_{L^2_x}^3
    -\dfrac{c}{2}\gamma\delta\Big).\label{2.23}
\end{eqnarray}
The last term of (\ref{2.23}) is negative if
\begin{eqnarray}
&N^{-\frac{3}{2}}\>\|Iu(t)\|_{L^2_x}, N^{-3}\>\|Iu(t)\|_{L^2_x}^2,
N^{-\frac{3}{2}-s}\>\|f\|_{H^s},\nonumber\\
&N^{-3-s}\>\|Iu(t)\|_{L^2_x}\>\|f\|_{H^s},
N^{-\beta}\delta^{-\frac{2}{3}+}\>\|Iu(t)\|_{L^2_x}^3= o(N).
\label{2.24}
\end{eqnarray}
Now we consider it by iteration. When $t=0$, by (\ref{II}) and
noting that $\delta\sim N^{3s-}$, (\ref{2.24}) is satisfied when
$s>-\dfrac{1}{2}$. We just check the last term of (\ref{2.24})
here, which is sufficient to show
$$
N^{-\beta}\cdot N^{-2s+}\cdot N^{-3s}\|u_0\|_{H^s}\lesssim N^{0-},
$$
that is exactly, $s>-\dfrac{1}{2}$. Therefore, (\ref{2.23}) yields
\begin{equation}\label{iteration}
E^4_I(t)\leq e^{-\gamma t}E^4_I(0) +K(f), \quad \mbox{for any }
t\in[0,\delta].
\end{equation}
By using (\ref{E3}), (\ref{E4}), (\ref{iteration}) and Lemma 3.2
again, we get
\begin{equation}\label{2.25}
E^2_I(t)
    \leq
    e^{-\gamma t}E^4_I(0)
    +K(f)+CN^{-\frac{3}{2}}\>\|Iu(t)\|_{L^2_x}^3
 +CN^{-3}\>\|Iu(t)\|_{L^2_x}^4
\end{equation}
for any $t\in [0,\delta]$. Note that
$$
E^4_I(0)\leq C\|Iu_0\|_{L^2_x}^2+
CN^{-\frac{3}{2}}\>\|Iu_0\|_{L^2_x}^3
 +CN^{-3}\>\|Iu_0\|_{L^2_x}^4,
$$
then by (\ref{2.25}), it is easy to see that for $t\in [0,\delta]$,
\begin{equation}\label{2.27}
E^2_I(t)
    \leq
    2e^{-\gamma\delta}E^4_I(0)+2K(f)
\end{equation}
by choosing $N$, dependent on $\|u_0\|_{H^s},\|f\|_{H^s}$, large
enough. By (\ref{2.27}), (\ref{2.24}) is true if we take $t=\delta$.
Repeating the process above, iterating (\ref{iteration}), we
conclude that (\ref{2.27}) holds for any $t>0$. Therefore, by
(\ref{II}) again, we have
\begin{equation}\label{GEu}
\limsup\limits_{t\rightarrow+\infty}\|u(t)\|_{H^s}
    \leq
    \sqrt{2K(f)}.
\end{equation}

We state our results in this subsection.
\begin{thm}
Let $\gamma>0, f\in \dot{H}^s(\T)$ and $s>-\dfrac{1}{2}$, then (\ref{KdV})
(\ref{initial data}) are globally well-posed in $\dot{H}^s(\T)$. Moreover,
the solution $u(t)$ satisfies (\ref{GEu}).
\end{thm}

By Theorem 3.4, (\ref{KdV}) (\ref{initial data}) define a continuous
semigroup $S(t)$ on $H^s(\T)$, such that $u(t)=S(t) u_0\in
C([0,+\infty);H^s(\T))$ is the solution with the initial function
$u_0$. Moreover, we have
\begin{cor}
Let $s>-\dfrac{3}{4}$, then the solution map $S(t)$ associated with
system (\ref{KdV}) possesses a bounded absorbing ball in $H^s(\T)$,
with the radius given by (\ref{GEu}).
\end{cor}

 \vspace{0.3cm}
\section{Split of the Solution}
For the existence of the global attractor, we need to show the
asymptotic compactness of the solution map in $H^s(\T)$. Since the
KdV equation possesses no conservation law for in $H^s(\T)$ for
$s<0$ and lacks sufficient regularity of the solutions, we apply the
asymptotic smoothing effect via a suitable decomposition of the
solution map to obtain it.

For this purpose, we fix a large number $N$ (which will be chosen
later and may be different from the one in Section 3) and split
the solution $u$ into two parts as
$$
u=v+w,
$$
where
\begin{eqnarray}
&&\partial_{t}v+\partial_{x}^{3}v+\dfrac{1}{2}\partial_xv^2+\gamma v
  =-\dfrac{1}{2}P_N(\partial_xu^2-\partial_xv^2)+f,\label{Eqv}\\
&&\partial_{t}w+\partial_{x}^{3}w+\gamma w
  =-\dfrac{1}{2}Q_N(\partial_xu^2-\partial_xv^2),\label{Eqw}\\
&&v(x,0)=P_Nu_0(x),\qquad w(x,0)=Q_Nu_0(x).\label{initial data vw}
\end{eqnarray}

The local well-posedness of the systems (\ref{Eqw})--(\ref{initial
data vw}) can be proved by employing the bilinear estimates
(\ref{L21}) and the standard process of the fixed point argument.
Especially, taking the initial time ``$t_0$" (under the assumption
of existence), we have the estimate
\begin{equation}\label{Ew}
\|w\|_{Y^s} \lesssim
\|w(t_0)\|_{H^s},
\end{equation}
with the lifetime $\delta\in (0,1)$ depending on $\gamma$,
$\|w(t_0)\|_{H^s}$ and $\|f\|_{H^s}$ but independent of $N$.
Further, by (\ref{Ew}), we have
\begin{equation}\label{Ew2}
\sup\limits_{t\in [t_0-\delta,t_0+\delta]} \|w(t)\|_{H^s} \lesssim
\|w(t_0)\|_{H^s}.
\end{equation}

For $v$, by (\ref{LSE}) (for $N=1$) and (\ref{Ew}), we have
\begin{eqnarray}
&\|v\|_{Y^s} \lesssim
\|v(t_0)\|_{H^s}+\|u(t_0)\|_{H^s}+\|f\|_{H^s},\label{Ev}\\
&\sup\limits_{t\in [t_0-\delta,t_0+\delta]} \|v(t)\|_{H^s} \lesssim
\|v(t_0)\|_{H^s}+\|u(t_0)\|_{H^s}+\|f\|_{H^s}.\label{Ev2}
\end{eqnarray}

\subsection{Decay of $w$ in $H^s(\T)$}
Noting that $w=Q_Nw$, we rewrite (\ref{Eqw}) into
\begin{equation}\label{3.6}
    \partial_{t}w+\partial_{x}^{3}w+\gamma w
  =Q_N[ww_x-(uw)_x],
\end{equation}
and drive the energy equation of (\ref{3.6}) in $H^s(\R)$ to find
\begin{eqnarray}
\|J^s_x w(t+\delta)\|_{L^2}^2 &=&e^{-\gamma\delta}\|J^s_x
w(t)\|_{L^2}^2 +2\displaystyle\int_0^\delta
e^{-\gamma(\delta-t')}\!\int_\T\!J^s_x[ww_x-(uw)_x]
\cdot J^s_x w\,dx dt'\nonumber\\
&&-\gamma \displaystyle\int_0^\delta e^{-\gamma(\delta-t')}\|J^s_x
w\|_{L^2}^2\,dt', \label{Ew3}
\end{eqnarray}
where we have omitted the variable $t'-\delta+t$ of the functions
inside the time integral for short. Since for each $f\in
X_{s,\frac{1}{2}}^{[t,t+\delta]}$, there exists an $\tilde{f}\in
X_{s,\frac{1}{2}}$ such that $f\big|_{t'\in
[t,t+\delta]}=\tilde{f}\big|_{t'\in [t,t+\delta]}$ and $
\|f\|_{X_{s,\frac{1}{2}}^{[t,t+\delta]}}
=\left\|\tilde{f}\right\|_{X_{s,\frac{1}{2}}}$, we may replace
$u,w$ by $\tilde{u},\tilde{w}$ in the following procedure. But we
remove the tilde $\tilde{}$ again for simplicity.

First, by (\ref{Ew2}) we have
\begin{equation}\label{3.10}
\gamma\displaystyle\int_0^\delta e^{-\gamma(\delta-t')} \|J^s_x
w\|_{L^2}^2\,dt' \geq c(1-e^{-\gamma\delta})\>\|w(t)\|_{H^s_x}^2.
\end{equation}
For the second term in (\ref{Ew3}), We write
$$
\displaystyle\int_0^\delta e^{-\gamma(\delta-t')}\!
\int_\T\!J^s_x[ww_x-(uw)_x]\cdot J^s_x w\,dx\,dt'=J_1+J_2,
$$
where
$$
\begin{array}{c}
J_1=\dfrac{1}{2}\displaystyle\int_0^\delta e^{-\gamma(\delta-t')}\!
\int_\T\!J^s_x\partial_x(w^2)\cdot J^s_x w\,dx\,dt',\\
J_2=-\displaystyle\int_0^\delta e^{-\gamma(\delta-t')}\!
\int_\T\!J^s_x\partial_x(uw)\cdot J^s_x w\,dx\,dt'.
\end{array}
$$
We rewrite $J_2$ again by
$$
J_2=J_{21}+J_{22}+J_{23}
$$
where
$$
\begin{array}{cll}
J_{21} &\!\!=&\!\! \dfrac{1}{2}\displaystyle\int_0^\delta
e^{-\gamma(\delta-t')}\!
\int_\T\!P_{N^\epsilon} u\cdot \partial_x(J^s_x w)^2\,dx\,dt',\\
J_{22} &\!\!=&\!\! -\!\displaystyle\int_0^\delta \!\!
e^{-\gamma(\delta-t')}\!\! \int_\T\!\partial_x[J^s_x,P_{N^\epsilon} u]w
\cdot J^s_x w\,dx\,dt',\\
J_{23} &\!\!=&\!\! -\displaystyle\int_0^\delta
e^{-\gamma(\delta-t')}\! \int_\T\!\partial_xJ^s_x(Q_{N^\epsilon} u\cdot
w)\cdot J^s_x w\,dx\,dt',
\end{array}
$$
where $0<\epsilon\ll 1$, $N\gg1$, and the commutator $[A,B]=AB-BA$.
First,
\begin{lem}
For any $s\geq -\dfrac{1}{2}$, the functions $z_0,z$ are zero
$x-$mean for all $t$, and $z_0\in X_{0,\frac{1}{2}}^\delta$,
$z\in X_{s,\frac{1}{2}}^\delta$ with $z=Q_Nz$, then,
\begin{eqnarray}
&\displaystyle\int_0^\delta e^{-\gamma(\delta-t')}\! \int_\T\!z_0\cdot
\partial_x(J^s_x z)^2\,dx\,dt' \lesssim
N^{-\frac{1}{2}}\>\|z_0\|_{X_{0,\frac{1}{2}}^\delta}\>
\|z\|_{X_{s,\frac{1}{2}}^\delta}^2.\label{3.12}
\end{eqnarray}
\end{lem}
\begin{proof}
By duality and (\ref{L21'}), the left-hand side of (\ref{3.12}) is
controlled by
$$
\|z_0\|_{X_{0,\frac{1}{2}}^\delta}\> \|\varphi \partial_x(J^s_x
z)^2\|_{X_{0,-\frac{1}{2}}^\delta} \lesssim
\|z_0\|_{X_{0,\frac{1}{2}}^\delta}\> \|J^s_x
z\|_{X_{0,\frac{1}{2}}^\delta}\> \|J^s_x
z\|_{X_{-\frac{1}{2},\frac{1}{2}}^\delta},
$$
where $\varphi= \chi_{[0,\delta]}e^{-\gamma(\delta-t')}$. Then the
result follows by noting $z=Q_Nz$. \hfill$\Box$
\end{proof}

By this lemma, (\ref{LSE}) (for $N=1$) and (\ref{Ew}), we have
\begin{equation}\label{J21}
J_{21}\lesssim
N^{-\frac{1}{2}-s\epsilon}\>\|u\|_{X_{s,\frac{1}{2}}^\delta}\>
\|w\|_{X_{s,\frac{1}{2}}^\delta}^2 \lesssim
N^{-\frac{1}{2}-s\epsilon}\>(\|u(t)\|_{H^s}+\|f\|_{H^s})\>
\|w(t)\|_{H^s}^2.
\end{equation}

\begin{lem}
For any $s\in \R$, the functions $z_0,z$ are zero
$x-$mean for all $t$, and $z_0\in X_{\frac{1}{2},\frac{1}{2}}^\delta$ and $z\in
X_{s,\frac{1}{2}}^\delta$ with $z_0=P_{\ll
N}z_0,z=Q_Nz$, then,
\begin{eqnarray}
&\displaystyle\int_0^\delta \!\! e^{-\gamma(\delta-t')}\!\!
\int_\T\!\partial_x[J^s_x,z_0]z \cdot J^s_x z\,dx\,dt' \lesssim
N^{-1}\|z_0\|_{X_{\frac{1}{2},\frac{1}{2}}^\delta}\>
\|z\|_{X_{s,\frac{1}{2}}^\delta}^2.\label{L4.2}
\end{eqnarray}
\end{lem}
\begin{proof}
It follows easily by duality and (\ref{L5.4'}), where we set
$\varphi= \chi_{[0,\delta]}e^{-\gamma(\delta-t')}$. \hfill$\Box$
\end{proof}

By this lemma, (\ref{LSE}) and (\ref{Ew}), we have
\begin{equation}\label{J22}
J_{22}\lesssim
N^{-1+(\frac{1}{2}-s)\epsilon}\>\|u\|_{X_{s,\frac{1}{2}}^\delta}\>
\|w\|_{X_{s,\frac{1}{2}}^\delta}^2 \lesssim
N^{-1+(\frac{1}{2}-s)\epsilon}\>(\|u(t)\|_{H^s}+\|f\|_{H^s})\>
\|w(t)\|_{H^s}^2.
\end{equation}

On the other hand, by the duality, (\ref{L21'}), (\ref{LSE}) and
(\ref{Ew2}), we have
\begin{eqnarray}
J_{23} &\lesssim& \|\varphi \partial_xJ^s_x(Q_{N^\epsilon}u\cdot
w)\|_{X_{0,-\frac{1}{2}}}\>
\|J^s_x w\|_{X_{0,\frac{1}{2}}}\nonumber\\
&\lesssim& \left(\|Q_{N^\epsilon}u\|_{X_{s,\frac{1}{2}}^\delta}
\|w\|_{X_{-\frac{1}{2},\frac{1}{2}}^\delta}
+\|Q_{N^\epsilon}u\|_{X_{-\frac{1}{2},\frac{1}{2}}^\delta}
\|w\|_{X_{s,\frac{1}{2}}^\delta}\right)\|w\|_{X_{s,\frac{1}{2}}^\delta}
\nonumber\\
&\lesssim&
N^{-(s+\frac{1}{2})\epsilon}\|Q_{N^\epsilon}u\|_{X_{s,\frac{1}{2}}^\delta}
\|w\|_{X_{s,\frac{1}{2}}^\delta}^2\nonumber\\
&\lesssim&
N^{-(s+\frac{1}{2})\epsilon}(\|u(t)\|_{H^s}+\|f\|_{H^s})
\|w(t)\|_{H^s}^2 .\label{J23}
\end{eqnarray}

Similar to $J_{23}$, we have
\begin{equation}\label{J1}
    J_1\lesssim
N^{-(s+\frac{1}{2})} \|w(t)\|_{H^s}^3 .
\end{equation}

Summing up (\ref{J21}), (\ref{J22})--(\ref{J1}), we have
\begin{equation}\label{3.8}
   \displaystyle\int_0^\delta e^{-\gamma(\delta-t')}\!
\int_\T\!J^s_x[ww_x-(uw)_x]\cdot J^s_x w\,dx\,dt' \lesssim
\epsilon(N)\>(1+\|w(t)\|_{H^s_x})\>\|w(t)\|_{H^s_x}^2
\end{equation}
for some $\epsilon(N)=o(N)$, where we have used (\ref{GEu}).

Inserting (\ref{3.10}), (\ref{3.8}) into (\ref{Ew3}), we have
\begin{eqnarray}
&&\|J^s_x w(t+\delta)\|_{L^2}^2\nonumber\\
&\leq& e^{-\gamma\delta}\|J^s_x
w(t)\|_{L^2}^2+\left[\epsilon(N)(1+\|w(t)\|_{H^s_x})-
c(1-e^{-\gamma\delta})\right]\|w(t)\|_{H^s_x}^2.\label{3.11}
\end{eqnarray}

Note that $\|w(0)\|_{H^s}=\|Q_Nu_0\|_{H^s}\leq \|u_0\|_{H^s}$,
therefore, we observe that the last term in (\ref{3.11}) is negative
when $t=0$ by choosing $N$, dependent only on
$\gamma,\delta,\|u_0\|_{H^s}$, large enough. Hence, for $t\in
[0,\delta]$,
\begin{equation}\label{Dyw}
\|w(t)\|_{H^s_x} \leq e^{-\gamma t}\| Q_Nu_0\|_{H^s}^2 .
\end{equation}
By iteration, we conclude that (\ref{Dyw}) holds for any $t>0$.
\begin{prop}
Let $\gamma>0$, $u_0\in \dot{H}^s(\T)$, $u$ is the solution of
(\ref{KdV}) (\ref{initial data}) given in Proposition 3.1, then
(\ref{Eqw}) (\ref{initial data vw}) are global well-posedness in
$\dot{H}^s(\T)$ and the solution satisfies the decay estimate
(\ref{Dyw}). Moreover, for any $T>0$,
\begin{equation}\label{Ew4}
\|w\|_{X_{s,\frac{1}{2}}^T} \leq C\|Q_N u_0\|_{H^s},
\end{equation}
for some constant $C$ independent of $T$.
\end{prop}
\begin{proof}
We only need to show (\ref{Ew4}), which follows from (\ref{Ew}),
(\ref{Dyw}) and the estimate
\begin{equation}\label{Omega0}
\|f\|_{X_{s,b}^{\Omega}}\leq \|f\|_{X_{s,b}^{\Omega_0}}+
\|f\|_{X_{s,b}^{\Omega/\Omega_0}}.
\end{equation}
for any $f\in X_{s,b}^{\Omega}$, $\Omega_0\subset \Omega$ (see
\cite{LW2} for the proof in the real-line case). \hfill$\Box$
\end{proof}

\subsection{Regularity of $v$}

We obtain the global well-posedness of $v$, by the global
well-posedness of $u$ and $w$. Moreover, by (\ref{GEu}) and
(\ref{Dyw}), we have the bound
\begin{equation}\label{GEv}
\limsup\limits_{t\rightarrow+\infty}\|v(t)\|_{H^s}
    \leq
    2K(f).
\end{equation}

Now we consider the regularity of $v$. For this purpose, we split
$v$ again into two parts by writing
$$
v=y+Q_Nv,
$$ where $y=P_Nv=P_Nu$. Then by (\ref{LSE}) and (\ref{GEu}), we have
\begin{equation}\label{3.14}
\|y\|_{L^\infty_{t}H^m_{x}},\|y\|_{X_{m,\frac{1}{2}}^\delta}\lesssim
N^{m-s}, \quad \mbox{ for any } m\geq s.
\end{equation}

Therefore, we only focus on the regularity of $Q_Nv$. Note that the
method used in \cite{Goubet} is seemly not suitable for this
situation, since we hardly give enough estimates on $y_t$ in the low
regularity case. For this purpose, we employ the one used in
\cite{M} and \cite{Ts2}. First, we introduce the functions
\begin{equation}\label{g}
\hat{g}(\xi)=\dfrac{\hat{f}}{i\xi^3+\gamma},\qquad g_N=Q_Ng,
\end{equation}
then $g$ is the solution of
$$
\partial_{x}^{3}g+\gamma g=f.
$$
Let $z=Q_Nv-g_N$, then
\begin{eqnarray}
&&\partial_{t}z+\partial_{x}^{3}z+\gamma z
  =-\dfrac{1}{2}\partial_xQ_N(v^2),\label{Eqz}\\
&& z(0)=-g_N.\label{initial data z}
\end{eqnarray}
Now we turn to prove that $z$ is uniformly bounded in $H^{s+3}(\T)$
when $f\in H^{s}(\T)$. For this, some lemmas are needed.
\begin{lem} Let $g_N$ defined in (\ref{g}) for $f\in H^s(\T)$
and $l<\dfrac{3}{2}+s$, then
\begin{equation}\label{EgN}
\|g_N\|_{H^{s+3}}\leq \|Q_Nf\|_{H^{s}};\qquad
\|g_N\|_{X_{l,\frac{1}{2}}^\delta}\lesssim\|Q_Nf\|_{H^{s}}.
\end{equation}
\end{lem}
\begin{proof}
The first term follows from the definition (\ref{g}). On the other hand,
\begin{eqnarray*}
\left\|\psi(\tau/\delta)g_N\right\|_{X_{l,\frac{1}{2}}}
&=&\delta\left\|\langle\tau-\xi^3\rangle^{\frac{1}{2}}
\langle\xi\rangle^{l}
\hat{\psi}(\delta\tau)\widehat{g_N}\right\|_{L^2_{\xi \tau}}\\
&\lesssim&
\delta\left\|\langle\tau\rangle^{\frac{1}{2}}
\langle\xi\rangle^{l}
\hat{\psi}(\delta\tau)\widehat{g_N}\right\|_{L^2_{\xi \tau}}
+\delta\left\|\langle\xi\rangle^{l+\frac{3}{2}}
\hat{\psi}(\delta\tau)\widehat{g_N}\right\|_{L^2_{\xi \tau}}\\
&=& \delta\left\|\langle\tau\rangle^{\frac{1}{2}}
\hat{\psi}(\delta\tau)\right\|_{L^2_{\tau}}
\left\|g_N\right\|_{H^{l}}
+\delta\left\|\hat{\psi}(\delta\tau)\right\|_{L^2_{\tau}}
\left\|g_N\right\|_{H^{l+\frac{3}{2}}}\\
&\lesssim&
\left\|\psi\right\|_{H^{\frac{1}{2}}}
\left\|g_N\right\|_{H^{s+3}}.
\end{eqnarray*}
Then the second term follows from the result of the first term.
\hfill$\Box$
\end{proof}

\begin{lem}
For $s\geq -\dfrac{1}{2}$, $f\in \dot{H}^s(\T)$, $\gamma\in \R$,
(\ref{Eqz}) (\ref{initial data z}) are local well-posedness in
$\dot{H}^s(\T)$ with some lifetime $\delta$ which depends on
$\gamma,\|u_0\|_{H^s},\|f\|_{H^s}$ but is independent of $N$.
Moreover,
\begin{equation}\label{Ez1}
\|z\|_{Y^s}\lesssim K_1,
\end{equation}
for some constant $K_1$ independent of $N$.
\end{lem}
\begin{proof}
It follows easily from (\ref{L21}), (\ref{Ev}), (\ref{GEu}),
(\ref{GEv}) and the standard fixed point argument. \hfill$\Box$
\end{proof}

Next we prove the regurality of $z$.
\begin{prop}
For $s>-\dfrac{1}{2}$, $f\in \dot{H}^s(\T)$, $\gamma\in \R$, the
solution $z$ obtained in Lemma 4.5 belongs to $Y^{s+3} \subset
C([-\delta,\delta]; H^{s+3}(\T))$ with some lifetime $\delta$,
which depends on $\gamma$, $\|u_0\|_{H^s}$, $\|f\|_{H^s}$ but is
independent of $N$. Moreover,
\begin{equation}\label{Ez2}
\|z\|_{Y^{s+3}}\lesssim
K_2(N)+\|z(0)\|_{H^{s+3}}.
\end{equation}
\end{prop}

To prove this proposition, we need the following results.

\begin{lem}
For any $r\geq 0$, and $h\in H^{r}(\T)$, the following bilinear
estimate holds,
$$
 \|\psi(t/\delta)\partial_x(h^2)\|_{Z^r}
 \leq
 C\|h\|_{H^r}^2.
$$
\end{lem}
\begin{proof}
By replacing $Z^s$ by $X_{s,-\frac{1}{2}+}$, then it is indeed a
consequence of Lemma 2.5 in \cite{Ts2}, although it is given in
the real-line case. \hfill$\Box$
\end{proof}

{\it \noindent Proof of Proposition 4.6.\quad} Rewrite the
nonlinearity as
\begin{equation}\label{v^2}
v^2= y^2+z^2+g_N^2+2yz+2yg_N+2zg_N,
\end{equation}
then for some $\mu>0$,
\begin{eqnarray}
\|\partial_x(v^2)\|_{Z^{s+3}} &\lesssim&
C\|g_N\|_{H^{s+3}}^2+\delta^\mu
\left(\|y\|_{X_{s+3,\frac{1}{2}}}+\|z\|_{X_{s+3,\frac{1}{2}}}
\right)\nonumber\\
&&
\cdot\left(\|y\|_{X_{s,\frac{1}{2}}}+\|z\|_{X_{s,\frac{1}{2}}}+\|g_N\|_{H^{s+3}}
\right)\nonumber\\
&\lesssim& C\|f\|_{H^s}^2+\delta^\mu
\left(\|y\|_{X_{s+3,\frac{1}{2}}}+\|z\|_{X_{s+3,\frac{1}{2}}}\right)
,\label{Ev^2}
\end{eqnarray}
where we used (\ref{L21'''}) to treat the first, second and fourth
terms in (\ref{v^2}), used Lemma 4.7 to treat the third term and
used (\ref{L25}) to treat the fifth, sixth terms. Taking the
$Y^{s+3}$-norm onto the two sides of the Duhamel's integral equation
of (\ref{Eqz}), then the results easily follow from (\ref{3.14}),
(\ref{EgN}) and (\ref{Ez1}). \hfill$\Box$

We drive the energy equation of $z$ in $H^{s+3}(\T)$ to find
\begin{eqnarray}
&\|z(t+\delta)\|_{H^{s+3}}^2=e^{-\gamma\delta}\|z(t)\|_{H^{s+3}}^2\nonumber\\
&+\displaystyle\int_0^\delta e^{-\gamma(\delta-t')}
\!\left\{-\!\int_\T\!\partial_x^4J^s_x(v^2) \cdot \partial_x^3J^s_x
z\,dx-\gamma \|z\|_{H^{s+3}}^2\right\}\,dt',\label{3.29}
\end{eqnarray}
where we have omitted the variable $t'-\delta+t$ of the functions
inside the time integral. In the following, we will prove that
\begin{equation}\label{3.30}
\displaystyle\int_0^\delta e^{-\gamma(\delta-t')}
\!\!\int_\T\!\partial_x^4J^s_x(v^2) \cdot \partial_x^3J^s_x z\,dx dt'
\lesssim K_3(N)+\epsilon(N)\|z(t)\|_{H^{s+3}}^2.
\end{equation}
For this, we need the following lemma.
\begin{lem}
For any $s>-\dfrac{1}{2}$,
and $z_0\in X_{s,\frac{1}{2}}^\delta$, $z\in X_{s+3,\frac{1}{2}}^\delta$ with
$z=Q_Nz$,
we have
\begin{eqnarray}
&\displaystyle\int_0^\delta e^{-\gamma(\delta-t')}
\!\!\int_\T\!\partial_x^4J^s_x(z_0\cdot z) \cdot \partial_x^3J^s_x
z\,dx dt' \lesssim
\epsilon(N)\>\|z_0\|_{X_{s,\frac{1}{2}}^\delta}\>
\|z\|_{X_{s+3,\frac{1}{2}}^\delta}^2\nonumber\\
&+ \epsilon(N)\>\|z\|_{X_{s,\frac{1}{2}}^\delta}\>
\|z_0\|_{X_{s+3,\frac{1}{2}}^\delta}\>
\|z\|_{X_{s+3,\frac{1}{2}}^\delta}.\label{4.35}
\end{eqnarray}
\end{lem}
\begin{proof}
Combining Lemma 4.1 and Lemma 4.2, we have
$$
\displaystyle\int_0^\delta e^{-\gamma(\delta-t')}
\!\!\int_\T\!\partial_x^4J^s_x(P_{N^\epsilon}
z_0\cdot z) \cdot \partial_x^3J^s_x
z\,dx dt'
\lesssim
\epsilon(N)\|z_0\|_{X_{s,\frac{1}{2}}^\delta}\>
\|z\|_{X_{s+3,\frac{1}{2}}^\delta}^2,
$$
where $\epsilon(N)=o(N)$. On the other hand, by the duality and (\ref{L21'}),
$$
\begin{array}{lll}
&&\displaystyle\int_0^\delta e^{-\gamma(\delta-t')}
\!\!\int_\T\!\partial_x^4J^s_x(Q_{N^\epsilon}
z_0\cdot z) \cdot \partial_x^3J^s_x
z\,dx dt'\\
&\lesssim& \|\varphi \partial_x^4J^s_x(Q_{N^\epsilon}z_0\cdot
z)\|_{X_{0,-\frac{1}{2}}}\>
\|\partial_x^3J^s_x z\|_{X_{0,\frac{1}{2}}}\\
&\lesssim& \left(\|Q_{N^\epsilon}z_0\|_{X_{s+3,\frac{1}{2}}^\delta}
\|z\|_{X_{-\frac{1}{2},\frac{1}{2}}^\delta}
+\|Q_{N^\epsilon}z_0\|_{X_{-\frac{1}{2},\frac{1}{2}}^\delta}
\|z\|_{X_{s+3,\frac{1}{2}}^\delta}\right)\|z\|_{X_{s+3,\frac{1}{2}}^\delta}
\nonumber\\
&\lesssim& \epsilon(N)\left(\|z_0\|_{X_{s+3,\frac{1}{2}}^\delta}
\|z\|_{X_{s,\frac{1}{2}}^\delta}
+\|z_0\|_{X_{s,\frac{1}{2}}^\delta}
\|z\|_{X_{s+3,\frac{1}{2}}^\delta}\right)\|z\|_{X_{s+3,\frac{1}{2}}^\delta}.
\end{array}
$$
Summing up the two estimates above, we have (\ref{4.35}).
\hfill$\Box$
\end{proof}

Now we return to the proof of (\ref{3.30}). By (\ref{v^2}), we split
$\displaystyle\int_0^\delta e^{-\gamma(\delta-t')}
\!\!\int_\T\!\partial_x^4J^s_x(v^2)\\
\cdot \partial_x^3J^s_x z\,dx dt'$ into four parts by writing
$$
\begin{array}{l}
I_1=\displaystyle\int_0^\delta e^{-\gamma(\delta-t')}
\!\!\int_\T\!\partial_x^4J^s_x(y^2)
\cdot \partial_x^3J^s_x z\,dx dt';\\
I_2=\displaystyle\int_0^\delta e^{-\gamma(\delta-t')}
\!\!\int_\T\!\partial_x^4J^s_x(g_N^2)
\cdot \partial_x^3J^s_x z\,dx dt';\\
I_3=2\displaystyle\int_0^\delta e^{-\gamma(\delta-t')}
\!\!\int_\T\!\partial_x^4J^s_x[(y+z)g_N]
\cdot \partial_x^3J^s_x z\,dx dt';\\
I_4=\displaystyle\int_0^\delta e^{-\gamma(\delta-t')}
\!\!\int_\T\!\partial_x^4J^s_x[(2y+z)z] \cdot \partial_x^3J^s_x z\,dx
dt'.
\end{array}
$$
Then we treat $I_1$ by duality, (\ref{L21'}), (\ref{3.14}), (\ref{Ez2})
and Young's inequality to find
$$
\begin{array}{ll}
I_1 &\lesssim
\|y\|_{X_{s,\frac{1}{2}}^\delta}\|y\|_{X_{s+3,\frac{1}{2}}^\delta}
\|z\|_{X_{s+3,\frac{1}{2}}^\delta}\\
&\lesssim N^3\left(K_2(N)+\|z(t)\|_{H^{s+3}}\right)\\
&\leq K(N)+\epsilon(N)\|z(t)\|_{H^{s+3}}^2.
\end{array}
$$
We treat $I_2$ by duality, Lemma 4.7, (\ref{EgN}), (\ref{Ez2}) and
Young's inequality to find
$$
I_2 \lesssim \|g_N\|_{H^{s+3}}^2\|z\|_{X_{s+3,\frac{1}{2}}^\delta}
\leq K(N)+\epsilon(N)\|z(t)\|_{H^{s+3}}^2.
$$
We treat $I_3$ by duality, (\ref{L25'}), (\ref{3.14}), (\ref{Ez2})
and Young's inequality to find
$$
\begin{array}{ll}
I_3&\lesssim
\left(\|y+z\|_{X_{-\frac{1}{4}+,\frac{1}{2}}^\delta}\|g_N\|_{H^{s+3}}+
\|y+z\|_{X_{s+3,\frac{1}{2}}^\delta}\|g_N\|_{H^{-\frac{1}{4}+}}\right)
\|z\|_{X_{s+3,\frac{1}{2}}^\delta}\\
&\lesssim
\|y\|_{X_{s+3,\frac{1}{2}}^\delta}\|g_N\|_{H^{s+3}}
\|z\|_{X_{s+3,\frac{1}{2}}^\delta}+
\epsilon(N)\|g_N\|_{H^{s+3}}
\|z\|_{X_{s+3,\frac{1}{2}}^\delta}^2\\
&\lesssim K(N)+\epsilon(N)\|z(t)\|_{H^{s+3}}^2.
\end{array}
$$
We treat $I_4$ by employing Lemma 4.8, (\ref{3.14}), (\ref{Ez1}),
(\ref{Ez2}) and Young's inequality to find
$$
\begin{array}{ll}
I_4\!\! &\lesssim \epsilon(N)
\left(\|y\|_{X_{s,\frac{1}{2}}^\delta}
\|z\|_{X_{s+3,\frac{1}{2}}^\delta}+
\|y\|_{X_{s+3,\frac{1}{2}}^\delta}
\|z\|_{X_{s,\frac{1}{2}}^\delta}+ \|z\|_{X_{s,\frac{1}{2}}^\delta}
\|z\|_{X_{s+3,\frac{1}{2}}^\delta}
\right)\|z\|_{X_{s+3,\frac{1}{2}}^\delta}\\
&\lesssim K(N)+\epsilon(N)\|z(t)\|_{H^{s+3}}^2.
\end{array}
$$
Then (\ref{3.30})
follows from the estimates above.

On the other hand, by (\ref{Ez2}) we have
$$
\gamma\displaystyle\int_0^\delta e^{-\gamma(\delta-t')}
\|z\|_{H^{s+3}}^2\,dt' \geq
c(1-e^{-\gamma\delta})\>\|z(t)\|_{H^{s+3}}^2-K_2(N).
$$
Therefore, combining this with (\ref{3.29}), (\ref{3.30}),  we
have
$$
\|z(t+\delta)\|_{H^{s+3}}^2\leq
e^{-\gamma\delta}\|z(t)\|_{H^{s+3}}^2 +K_3(N)
+\left[C\epsilon(N)-c(1-e^{-\gamma\delta})\right]\|z\|_{H^{s+3}}^2.
$$
Note that the last term is always negative if we choose $N$, dependent on
$s,\gamma, \delta,$ large
enough. So by iteration we have,
$$
\|z(t)\|_{H^{s+3}}^2\leq e^{-\gamma t}\|g_N\|_{H^{s+3}}^2 +K_3(N).
$$
Since $v=y+z+g_N$, combining it with (\ref{3.14}), (\ref{EgN}), we
have
\begin{equation}\label{Ev3}
    \|v(t)\|_{H^{s+3}}^2\leq e^{-\gamma t}\|Q_Nf\|_{H^{s}}^2
+K_4(N).
\end{equation}
\begin{prop}
For $s>-\dfrac{1}{2}$, $f\in H^s(\T)$, $\gamma>0$, (\ref{Eqv})
(\ref{initial data vw}) are global well-posedness in
$H^{s+3}(\T)$. Moreover,  the solution satisfies (\ref{Ev3}) and
\begin{equation}\label{Ev4}
\|v\|_{X_{l,\frac{1}{2}}^T}\lesssim K_5(T),
\end{equation}
for any $l<\dfrac{3}{2}+s$, $T>0$, where $K_5$ is dependent on
$s,\gamma, \|u_0\|_{H^s},\|f\|_{H^s},T.$
\end{prop}
\begin{proof}
We only need to see (\ref{Ev4}), which follows from (\ref{3.14}),
(\ref{EgN}), (\ref{Ez2}) and (\ref{Omega0}). \hfill$\Box$
\end{proof}

 \vspace{0.3cm}

\section{Existence of Global Attractor and Asymptotic \\ Smoothing Effect}

\subsection{Existence of Global Attractor in $\dot{H}^s(\T)$}
Based on Corollary 3.5, we need to show the asymptotically compact
of the solution map $S(t)$ in $H^s(\T)$ to prove the existence of
global attractor. Let $\{u_{0_n}\}_n$ be a bounded sequence of
initial data in $\dot{H}^{s}(\T)$ and the time sequence $\{t_n\}_n$
tending to infinity. Let $u_n(t)=S(t)u_{0_n}$ be the corresponding
solution of (\ref{KdV}) (\ref{initial data}) and write
$u_n(t)=v_n(t)+w_n(t)$, where $v_n(t),w_n(t)$ are the solutions of
(\ref{Eqv})--(\ref{initial data vw}) corresponding to the initial
condition $u_{0_n}$. The plan now is to show that $u_n(t_n)$ is
precompact in $H^s(\T)$, with $w_n(t_n)$ decay to zero and
$v_n(t_n)$ bounded in $H^{s+3}(\T)$ and precompact in $H^s(\T)$.

By Proposition 4.9, we first get that for any $T>0$,
\begin{equation}\label{51}
\left\{v_n(t_n+\cdot)\right\}_n \quad \mbox{ is bounded in }
C([0,T];H^{s+3}(\T)).
\end{equation}

We recall that, a sequence $\{f_n(t)\}_n$, for $t\in \Omega$, is
{\it uniformly equicontinuous} in a Banach space $X$, if for any
$\epsilon>0$, there exist an $\eta>0$, such that, for any $n\in \N$,
$t,t'\in \Omega$,
$$
\|f_n(t)-f(t')\|_X\leq \epsilon, \mbox{\quad if \quad} |t-t'|\leq
\eta.
$$

\begin{lem}
For $t\in [0,T]$, $\left\{v_n(t)\right\}_n$ is {\it uniformly
equicontinuous} in $H^{s}(\T)$.
\end{lem}
\begin{proof}
For any $ t,\eta\in \R$, we have
\begin{eqnarray}
\|v_n(t+\eta)-v_n(t)\|_{H^{s}_x} \leq
\displaystyle\left\|\langle\xi\rangle^{s}\int\!|e^{i\eta\tau}-1|
|\widehat{v_n}(\xi,\tau)|\,d\tau\right\|_{L^2_{\xi}} .\label{52}
\end{eqnarray}

Given $\epsilon>0$, since
$\|\langle\xi\rangle^{s}\widehat{v_n}\|_{L^2_\xi L^1_\tau}$ is
uniformly bounded, there exists a large
$\tau_0=\tau_0(\epsilon)>0$, such that
$$
\displaystyle\left\|\langle\xi\rangle^{s} \int_{|\tau|\geq
\tau_0}\! |\widehat{v_n}(\xi,\tau)|\,d\tau\right\|_{L^2_{\xi}}
\leq \epsilon.
$$
On the other hand,
$$
\sup\limits_{|\tau|\leq \tau_0}
\left|e^{i\eta\tau}-1\right|\lesssim |\eta\tau_0|,
$$
which leads to
$$
\displaystyle\left\|\langle\xi\rangle^{s} \int_{|\tau|\leq
\tau_0}\!
|e^{i\eta\tau}-1||\widehat{v_n}(\xi,\tau)|\,d\tau\right\|_{L^2_{\xi}}
\lesssim
|\eta\tau_0|\|\langle\xi\rangle^{s}\widehat{v_n}\|_{L^2_\xi
L^1_\tau} \leq \epsilon
$$
by choosing $\eta_0=\eta_0(\epsilon)$ small enough. Then the claimed
result follows from (\ref{52}). \hfill$\Box$
\end{proof}

Combining (\ref{51}) and Lemma 5.1, we come to the conclusion, by
Arzela-Ascoli's theorem, that there exists a function $\bar{u}$
such that for any $l<s+3$,
\begin{eqnarray}
\bar{u}\in C([0,T];H^{l}(\T))\cap C_w([0,T];H^{s+3}(\T)) \cap
L^\infty([0,T];H^{s+3}(\T)),\label{4.1}
\end{eqnarray}
and there exists a subsequence of $\{n\}$ (we also denote it by
$\{n\}$) such that for any $t\in \R$,
\begin{eqnarray}\label{4.2}
v_n(t_n+t)&\rightharpoonup& \bar{u}(t)\quad \mbox{ weakly in }
H^{s+3}(\T),\nonumber\\
&\rightarrow& \bar{u}(t)\quad \mbox{ strongly in }
H^{l}(\T) \quad\mbox{ for any } l<s+3\label{54}.
\end{eqnarray}
Moreover, by (\ref{Dyw}), we have for any $t\in \R, t_n\rightarrow +\infty$,
\begin{equation}\label{4.3}
w_n(t_n+t)\rightarrow 0\quad \mbox{ strongly in } H^{s}(\T) .
\end{equation}
Therefore,
$$
u_n(t_n+t)\rightharpoonup \bar{u}(t) \quad \mbox{ strongly in }
H^{s}(\T).
$$
 Hence, we establish the following result.
\begin{prop}
Let $\gamma>0$, $f\in \dot{H}^s(\T)$, the weakly damped, forced KdV
equation (\ref{KdV}) possesses a global attractor $\mathscr{A}$ in
$\dot{H}^s(\T)$, which is bounded in $H^{s+3}(\T)$ and compact in
$H^{l}(\T)$ for any $l<s+3$.
\end{prop}

\subsection{Compactness of the Global Attractor in $H^{s+3}(\T)$}

In this subsection, we prove that the attractor is in fact compact
in $H^{s+3}(\T)$. For this purpose, we just restrict the flow on the
global attractor and assume that the sequence of the initial data
$\{{u_0}_n\}_n$ belongs to $\mathscr{A}$. Since
$S(t)\mathscr{A}=\mathscr{A}$ for any $t\geq 0$, it is easy to see
that the corresponding trajectories $u_n(t)$ are uniformly bounded
in $H^{s+3}(\T)$.

We consider $u_n'(t)=\dfrac{d}{dt}u_n(t)$ in $H^s(\T)$. Note that by
(\ref{KdV}), $u_n'(t)$ are uniformly bounded in $H^s(\T)$, and
satisfy that
\begin{eqnarray}
u_n'(t_n+\cdot)&\rightarrow& \bar{u}'(\cdot),\quad \mbox{ in}\quad
C([0,T];H^{s'}(\T)) \quad\mbox{for any } s'<s;\nonumber\\
&\rightharpoonup& \bar{u}'(\cdot),\quad \mbox{ weakly in}\quad
C_w([0,T];H^{s}(\T)).\label{4.8}
\end{eqnarray}
Moreover, they satisfy the equation
\begin{equation}\label{Equ'}
  \partial_{t}u'+\partial_{x}^{3}u'+\gamma
  u'+\partial_x(uu')
  =0.
\end{equation}
By the fixed point argument process and the estimate (\ref{L21})
and (\ref{LSE}), we see that (\ref{Equ'}) is locally well-posed in
$Y^s$ for $s\geq-\dfrac{1}{2}$, and restricted on the time
interval $[-\delta,\delta]$,
\begin{equation}\label{Eu'}
    \|u'\|_{Y^{s}}
    \lesssim
    \|u'(0)\|_{H^{s}}.
\end{equation}
Especially, we have by the continuity and (\ref{4.8}) that for any
$t\geq 0$,
\begin{eqnarray}\label{4.11}
\|u_n'(t_n+t+\cdot)-\bar{u}'(t+\cdot)\|_{X_{-\frac{1}{2},\frac{1}{2}}^\delta}
\lesssim \|u_n'(t_n+t)-\bar{u}'(t)\|_{H^{-\frac{1}{2}}_x}
\rightarrow 0,
\end{eqnarray}
when $ t_n\rightarrow+\infty.$

Now we drive the energy equations of $u_n'(t)$ to obtain that
\begin{eqnarray}
   \|u_{n}'(t_{n})\|_{H^s}^2
   =
   e^{-2\gamma T}\|u_{n}'(t_{n}-T)\|_{H^s}^2
   -2\displaystyle\int_0^T \! \!
   e^{-2\gamma (T-t')}\int_\T
   \partial_{x}J^s_x(u_{n}u_{n}')\cdot
   J^s_xu_{n}'\, dx dt'.\label{4.12}
\end{eqnarray}
We plan to prove that
\begin{equation}\label{4.13}
\displaystyle\int_0^T \! \!
   e^{-2\gamma (T-t')}\int_\T
   \partial_{x}J^s_x(u_{n}u_{n}')\cdot
   J^s_xu_{n}'\, dx dt'
\rightarrow \displaystyle\int_0^T \! \!
   e^{-2\gamma (T-t')}\int_\T
   \partial_{x}J^s_x(\bar{u}\bar{u}')\cdot
   J^s_x\bar{u}'\, dx dt'.
\end{equation}
First, we obtain that when $ t_n\rightarrow+\infty$,
\begin{eqnarray}
\displaystyle\int_0^\delta \! \!
   e^{-2\gamma (T-t')}\int_\T
   \left[\partial_{x}J^s_x(u_{n}u_{n}')-\partial_{x}J^s_x(\bar{u}u_{n}')\right]\cdot
   J^s_xu_{n}'\, dx dt'
&\lesssim& \|u_{n}-\bar{u}\|_{X_{s,\frac{1}{2}}^\delta}
\|u_{n}'\|_{X_{s,\frac{1}{2}}^\delta}^2\nonumber\\
&\rightarrow&0,\label{4.14}
\end{eqnarray}
by duality, (\ref{L21'}), (\ref{4.2}) and the continuity given in Proposition
3.1.

On the other hand, applying the argument in Section 4.1, we write
$$
\displaystyle\int_0^\delta \! \!
   e^{-2\gamma (T-t')}\int_\T
   \partial_{x}J^s_x(\bar{u}u_{n}')\cdot
   J^s_xu_{n}'\, dx dt'\\
=H_1(\bar{u},u_{n}')+H_2(\bar{u},u_{n}'),
$$
where
$$
\begin{array}{l}
H_1(u,v)=-\dfrac{1}{2}\displaystyle\int_0^\delta \! \!
   e^{-2\gamma (T-t')}\int_\T
   u\cdot \partial_{x}(J^s_xv)^2
   \, dx dt';\\
H_2(u,v)=\displaystyle\int_0^\delta \! \!
   e^{-2\gamma (T-t')}\int_\T
   \partial_{x}[J^s_x,u]v\cdot
   J^s_xv\, dx dt'.
\end{array}
$$
Then by (\ref{L21'}) and (\ref{4.11}),  we have
\begin{eqnarray}\label{4.15}
&&H_1(\bar{u},u_{n}')-H_1(\bar{u},\bar{u}')\nonumber\\
&\lesssim& \|\bar{u}\|_{X_{\frac{1}{2}+s,\frac{1}{2}}^\delta}\>
\left\|\varphi\>\partial_{x}\left[J^s_x(u_{n}'-\bar{u}')\cdot
J^s_x(u_{n}'+\bar{u}')\right]
\right\|_{X_{-\frac{1}{2}-s,\frac{1}{2}}^\delta}\>
\nonumber\\
&\lesssim& \|\bar{u}\|_{X_{\frac{1}{2}+s,\frac{1}{2}}^\delta}\>
\|u_{n}'-\bar{u}'\|_{X_{-\frac{1}{2},\frac{1}{2}}^\delta}\>
\|u_{n}'+\bar{u}'\|_{X_{-\frac{1}{2},\frac{1}{2}}^\delta}
\nonumber\\
&\rightarrow &0,
\end{eqnarray}
where we note that $\|\bar{u}\|_{X_{l,\frac{1}{2}}^\delta}$ is
bounded for any $l<\dfrac{3}{2}+s$,
the proof can be given by the similar argument in Section
4.2.

We write $H_2(\bar{u},u_{n}')-H_2(\bar{u},\bar{u}')$ by
$$
\begin{array}{c}
\displaystyle\int_0^\delta \! \!
   e^{-2\gamma (T-t')}\int_\T
   \partial_{x}[J^s_x,\bar{u}](u_{n}'-\bar{u}')\cdot
   J^s_xu_{n}'\, dx dt'\\
+ \displaystyle\int_0^\delta \! \!
   e^{-2\gamma (T-t')}\int_\T
   \partial_{x}[J^s_x,\bar{u}] \bar{u}'\cdot
   J^s_x(u_{n}'-\bar{u}')\, dx dt'.
   \end{array}
$$
Then by duality, (\ref{L5.4''}) and (\ref{4.11}), we have
\begin{eqnarray}
&&
H_2(\bar{u},u_{n}')-H_2(\bar{u},\bar{u}')\nonumber\\
&\lesssim&
\left\|\varphi\>\partial_{x}[J^s_x,\bar{u}](u_{n}'-\bar{u}')
\right\|_{X_{\frac{1}{2}+s,-\frac{1}{2}}^\delta}
\|J^s_xu_{n}'\|_{X_{-\frac{1}{2}-s,\frac{1}{2}}^\delta}\nonumber\\
&&+ \left\|\varphi\>\partial_{x}[J^s_x,\bar{u}] \bar{u}'
\right\|_{X_{\frac{1}{2}+s,-\frac{1}{2}}^\delta}
\|J^s_x(u_{n}'-\bar{u}')\|_{X_{-\frac{1}{2}-s,\frac{1}{2}}^\delta}\nonumber\\
&\lesssim& \|\bar{u}\|_{X_{\frac{1}{2}+2s,\frac{1}{2}}^\delta}
\left(\|\bar{u}'\|_{X_{-\frac{1}{2},\frac{1}{2}}^\delta}
+\|u_{n}'\|_{X_{-\frac{1}{2},\frac{1}{2}}^\delta}\right)
\|u_{n}'-\bar{u}'\|_{X_{-\frac{1}{2},\frac{1}{2}}^\delta}
\nonumber\\
&\rightarrow &0.\label{4.16}
\end{eqnarray}

Note that we only restrict the integral domain with respect to time
on $[0,\delta]$, but it can be easily extended to $[0,T]$ by
dividing it into small intervals.

By (\ref{4.14})--(\ref{4.16}), we have (\ref{4.13}). By  the energy
equations of $\bar{u}'(t)$ in (\ref{4.12}) and let
$T\rightarrow+\infty$,
 we obtain that
$$
\limsup\limits_{t\rightarrow+\infty} \|u_{n}'(t_{n})\|_{H^s}^2 \leq
\|\bar{u}'(0)\|_{H^s}^2.
$$
Therefore,
$$
u_n'(t_n)\rightarrow \bar{u}'(0),\quad \mbox{ in }\quad H^{s}(\T)
$$
by combining (\ref{4.8}). This implies
$$
u_n(t_n)\rightarrow \bar{u}(0),\quad \mbox{ in }\quad H^{s+3}(\T).
$$
So we give the claim result in this subsection and thus finish the
proof of Theorem 1.1.


\vspace{0.3cm}
\begin{thebibliography}{99}


\bibitem{Ball}  J.\ M.\ Ball,
{\it Global attractors for damped semilinear wave equations. Partial
differential equations and applications}, Discrete and Contin.\
Dyn.\ Syst., {\bf 10}, 31--52, (2004).

\bibitem{B}  J.\ Bourgain,
{\it Fourier transform restriction phenomena for certain lattice
subsets and applications to nonlinear evolution equations. II. the
KdV equation}, Geom.\ Funct.\ Anal. {\bf 3}, 107--156, 209--262,
(1993)

\bibitem{CCT} M.\ Christ, J.\ Colliander and T.\ Tao,
{\it Asymptotics, frequency modulation, and low regularity
ill-posedness for canonical defocusing equations}, Amer.\ J.\ Math.,
{\bf 125} (6), 1235--1293, (2003).

\bibitem{CKSTT2} J.\ Colliander, M.\ Keel, G.\ Staffilani, H.\ Takaoka
and T.\ Tao, {\it Sharp global well-posedness for KdV and modified
Kdv on $\R$ and $\T$}, J.\ Amer.\ Math.\ Soc., {\bf 16}, 705--749,
(2003).


\bibitem{Gh2} J.\ M.\ Ghidaglia,
{\it Weakly damped forced Korteweg-de Vries equations behave as a
finite-dimensional dynamical system in the long time}, J. Differ.\
Eqns., {\bf 74}, 369--390, (1988).

\bibitem{Gh} J.\ M.\ Ghidaglia,
{\it A note on the strong convergence towards attractors of damped
forced KdV equations}, J. Differ.\ Eqns., {\bf 110}, 356--359,
(1994).

\bibitem{Goubet} O.\ Goubet,
{\it Asymptotic smoothing effect for weakly damped Korteweg-de Vries
equations},  Discrete and Contin.\ Dyn.\ Syst., {\bf 6}, 625--644,
(2000).

\bibitem{GR} O.\ Goubet and R.\ Rosa,
{\it Asymptotic smoothing and the global attractor of a weakly
damped KdV equation on the real line}, J.\ Differ.\ Eqns., {\bf 53},
25--53, (2002).

\bibitem{GH2} Guo, B.,  Huo, Z.:
{\it  The global attractor of the damped, forced generalized
Korteweg de Vries--Benjamin-Ono equaion in $L^2$}. Discrete and
Continuous Dynamical Systems {\bf 16} (1), 121--136, (2006).

\bibitem{KPV1} C.\ E.\ Kenig, G.\ Ponce and L.\ Vega,
{\it  A bilinear estimate with applications to the KdV equation},
J.\ Amer.\ Math.\ Soc., {\bf 9} (2), 573--603, (1996).

\bibitem{L} P.\ Laurencot,
{\it  Compact attractor for weakly damped driven Korteweg-de Vries
equations on the real line}, Czechoslovak Math.\ J., {\bf 48}
(1998), 85--94.

\bibitem{LW2} Y.\ Wu,
{\it  The Cauchy problem of the Schr\"{o}dinger-Korteweg-de Vries
system}, preprint.

\bibitem{M} L.\ Molinet,
{\it Global attractor and asymptotic smoothing effects for the
weakly damped cubic Schr\"{o}dinger equation in $L^2(\T)$}, Dynamics
of PDE, {\bf 6}, 15-34, (2009).

\bibitem{MoRo} I.\ Moise and R.\ Rosa,
{\it On the regularity of global attractor of a weakly damped,
forced Korteweg-de Vries equation}, Adv.\ Diff.\ Eqns., {\bf 2},
257--296, (1997).

\bibitem{MRW} I.\ Moise, R.\ Rosa and X.\ Wang,
{\it Attractors for non-compact semigroups via energy equations},
Nonlinearity, {\bf 11}, 1369--1393, (1998).

\bibitem{R} R.\ Rosa,
{\it The global attractor of a weakly damped, forced Korteweg-de
Vries equation in $H^1(\R)$}, Mat.\ Contemp., {\bf 19}, 129--152,
(2000).

\bibitem{Ts} K.\ Tsugawa,
{\it Existence of the global attractor for weakly damped, forced KdV
equation on Sobolev spaces of negative index}, Comm.\ Pure Appl.\
Anal., {\bf 3} (2), 513--528, (2004).

\bibitem{Ts2} K.\ Tsugawa,
{\it Global well-posedness for the KdV equations on the real line
with low regularity forcing terms}, Comm.\ Contemp.\ Math.\, {\bf 8}
(5), 681--713, (2006).




\end{thebibliography}
\end{document}